\documentclass[11pt,a4paper,reqno]{amsart}

\usepackage{mathpazo}
\usepackage{mdframed}
\usepackage[shortlabels]{enumitem}
\usepackage{tikz-cd}
\usepackage{appendix}
\usepackage{graphicx}
\usepackage{caption}
\usepackage{amsthm,amsmath,tikz}
\usepackage{mathrsfs}
\usepackage{xspace}
\usepackage[stable,perpage]{footmisc}
\usepackage{marginnote}
\usepackage{quiver}
\usepackage{amsfonts,amsmath,amssymb}
\usepackage{fancyhdr}
\usepackage{tikz-cd}
\usetikzlibrary{patterns,decorations.pathmorphing,shapes}
\usepackage{mathtools,stmaryrd}
\usetikzlibrary{decorations}
\usepackage{pgfplots}
\usepackage{bm}
\usepackage{hyperref,setspace,microtype}
\usepackage[margin=1in]{geometry}
\usepackage{amssymb} 
\def\acts{\curvearrowright}
\tikzset{
  symbol/.style={
    draw=none,
    every to/.append style={
      edge node={node [sloped, allow upside down, auto=false]{$#1$}}}
  }
}

\usepackage{mdframed}
\usepackage{setspace}

\theoremstyle{plain}
\newtheorem{thm}{Theorem}[section]
\newtheorem{lem}[thm]{Lemma}
\newtheorem{prop}[thm]{Proposition}
\newtheorem{cor}[thm]{Corollary}

\theoremstyle{definition}
\newtheorem{defn}[thm]{Definition}

\newtheorem*{cl}{Claim}

\newtheorem*{thm*}{Theorem}
\newtheorem*{defn*}{Definition}

\theoremstyle{remark}
\newtheorem*{pf}{Proof}
\newtheorem{rem}[thm]{Remark}

\newcommand{\R}{\mathbb{R}}

\newcommand{\Q}{\mathbb{Q}}

\newcommand{\OO}{\mathcal{O}}

\newcommand{\Z}{\mathbb{Z}}

\newcommand{\G}{\mathbb{G}}

\newcommand{\Spec}{\text{Spec}}

\newcommand{\rigidify}{%
  \mathrel{\begin{tikzpicture}[baseline=0ex, scale=0.15]
    \draw[line width=0.4pt] (0,0) -- (1,2);
    \draw[line width=0.4pt] (0.6,0) -- (1.6,2);
    \draw[line width=0.4pt] (1,2) -- (1.6,2);
    \draw[line width=0.4pt] (0,0) -- (0.6,0);
  \end{tikzpicture}}%
}

\newcommand{\Ext}{\text{Ext}}
\newcommand{\sym}{\text{Sym}}

\renewcommand{\hom}{\text{Hom}}

\newcommand{\PP}{\mathbb{P}}

\newcommand{\Def}{\text{Def}}

\newcommand{\QP}{\mathcal{Q}\mathcal{P}}

\newcommand{\hilb}{\mathcal{H}\text{ilb}}
\newcommand{\Ln}{L_{\underline{n}}}
\newcommand{\Lm}{L_{\underline{m}}}
\newcommand{\Lmn}{L_{\underline{m}+\underline{n}}}

\usepackage{stackengine}
\stackMath
\newcommand\tsup[2][2]{%
 \def\useanchorwidth{T}%
  \ifnum#1>1%
    \stackon[-.5pt]{\tsup[\numexpr#1-1\relax]{#2}}{\scriptscriptstyle\sim}%
  \else%
    \stackon[.5pt]{#2}{\scriptscriptstyle\sim}%
  \fi%
}

\usepackage{afterpage}

\title{Moduli of Multi-Uniformized Stacks and Seifert $\mathbb{G}_m^d$-Bundles}
\author{Zhengkai Pan}
\address{Department of Mathematics, Harvard University, Cambridge MA 02138, USA}
\email{zpan@math.harvard.edu}

\date{\today}
\begin{document}
\maketitle
\begin{abstract}
We introduce multi-uniformized stacks as a generalization of the Abramovich--Hassett construction of uniformized twisted varieties in \cite{AH} (hence also generalizing \cite{Hacking}). We prove the equivalence between the category of multi $\Q$-line bundles satisfying an analogue of  Koll\'ar's condition and the category of multi-uniformized twisted varieties, and we construct the corresponding moduli space. We then further broaden the framework to encompass Koll\'ar’s Seifert $\mathbb{G}_{m}^{d}$-bundles, showing that their moduli likewise coincide with those of $d$-uniformized $d$-cyclotomic orbispaces.
\end{abstract}

\setcounter{tocdepth}{1}
\tableofcontents
\section{Introduction}
Since the seminal introduction of the moduli space of stable curves $\overline{\mathcal{M}}_g$ by Deligne and Mumford in \cite{dm}, the pursuit of understandings of moduli spaces of varieties subject to prescribed numerical conditions, and of their natural compactifications, has become a central theme in algebraic geometry. In \cite{ksb88}, Koll\'ar and Shepherd--Barron initiated the study of moduli of stable surface pairs $(X, D)$ as a compactification of the moduli of surfaces of general type, thereby providing the first systematic higher‑dimensional analogue. In both the curve and surface settings, the relevant stable objects are required to carry semi‑log‑canonical (slc) singularities; in particular, they are reduced and satisfy Serre’s S2 condition.

In order to understand the moduli of surfaces when the boundary divisor $D$ is absent, an important condition on a stable object $X$ is that it has to be $\Q$-Gorenstein (i.e. the canonical divisor $K_X$ is a $\Q$-Cartier Weil divisor). If one wishes to study the deformation space of $X$, they could of course look at the space of all abstract deformations; but a more ``meaningful'' deformation space would be the tangent space to the moduli, so that the deformed surface still lies inside our moduli of interest (or continues to satisfy the stability conditions). These are $\Q$-Gorenstein deformations, or more commonly known as the KSB deformations. We denote corresponding deformation functor by
\[\text{Def}_{\text{KSB}}(X)\]
\subsection{Motivation}
This project was motivated by Paul Hacking's work \cite{Hacking}, where he constructed the canonical covering scheme on an slc surface $X$
\[\mathcal{Z} = \underline{\Spec}_X\bigoplus_{m=0}^{N-1}\omega_X^{[m]}\to X\]
(here $N$ is the Cartier index of $\omega_X$ at the $\Q$-Gorenstein singular point) and the canonical covering stack (also known as index 1 cover)
\[\mathcal{X} = \left[ \underline{\Spec}_X\bigoplus_{m=0}^{N-1}\omega_X^{[m]}\middle/ \mu_N\right]\xrightarrow{\pi} X\]
where $\pi$ is just the coarse space map, and here $\mathcal{X}$ is clearly a tame stack. The benefit of such a construction is that the KSB deformation space of $X$ has been proven in Proposition 3.7 of \cite{Hacking} to be equal to the abstract deformation space of the stack $\mathcal{X}$:
\[\text{Def}_{\text{KSB}}(X) =\text{Def}(\mathcal{X})\]
One way of understanding this equality is that deforming the stack $\mathcal{X}$ remembers the action, so it's more ``restrictive''. We remark here that on the deformed coarse space family $X'\to B$, the Koll\'ar's condition holds for $\omega_{X'/B}$, i.e.
\[\text{the formation of }\omega_{X'/B}^{[n]}\text{ for every $n\in\Z$ commutes with arbitrary base change}\]
Now, with this equivalence, we know that the KSB deformations of $X$ are controlled by the cotangent complex $\mathbb{L}_{\mathcal{X}/A}$ of $\mathcal{X}$, i.e.
\[\Def_{\text{KSB}}(X,M) = \Ext^1 (\mathbb{L}_{\mathcal{X}/A},\OO_{\mathcal{X}}\otimes M)\]
where $X$ sits over $\Spec A$ and $M$ is an $A$-module we extend along.

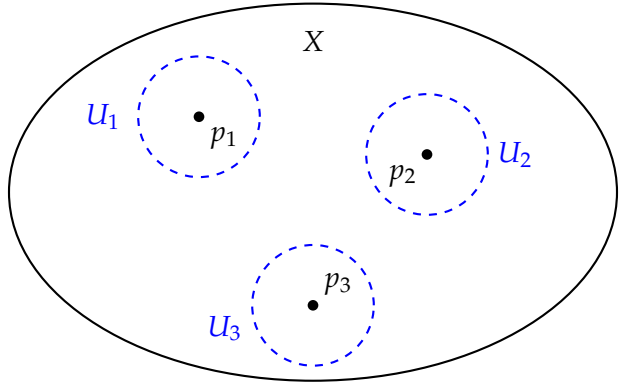
\begin{figure}
    \centering
\begin{tikzpicture}
  \draw[thick] (0,0) ellipse (4cm and 2.5cm);
  \node at (0,2) {\( X \)};
  
  \fill (-1.5,1) circle (2pt);
  \node[below right] at (-1.5,1) {\( p_1 \)};
  \draw[blue, dashed, thick] (-1.5,1) circle (0.8cm);
  \node[blue,left] at (-2.4,1) {\( U_1 \)};
  
  \fill (1.5,0.5) circle (2pt);
  \node[below left] at (1.5,0.5) {\( p_2 \)};
  \draw[blue, dashed, thick] (1.5,0.5) circle (0.8cm);
  \node[blue, right] at (2.3,0.5) {\( U_2 \)};
  
  \fill (0,-1.5) circle (2pt);
  \node[above right] at (0,-1.5) {\( p_3 \)};
  \draw[blue, dashed, thick] (0,-1.5) circle (0.8cm);
  \node[blue, below left] at (-0.8,-1.5) {\( U_3 \)};
  
\end{tikzpicture}
        \caption{Canonical covering with multiple singular points}
    \label{fig:1}
\end{figure}

Hacking's construction is local, and at each singular point $p\in X$ the construction depends on the Cartier index of $\omega_X$ at $p$. In particular, there might be multiple isolated singular points on the surface $X$. For example, in the cartoon Figure \ref{fig:1} above, $p_1,\,p_2,\,p_3$ are three isolated $\Q$-Gorenstein singular points. Hacking's construction of canonical cover in this case would be to take open neighborhoods $U_i$ around each $p_i$, construct the canonical covering stack $\mathcal{U}_i$ on each $U_i$ based on the Cartier index $N_i$ at $p_i$, and then glue them back together to get $\mathcal{X}$.

A global construction would be preferred because keeping track of the numbers $N_i$ for each singular point is something we wish to avoid when constructing the moduli space. This is exactly the work of Abramovich and Hassett in \cite{AH}. Starting from an arbitrary reflexive rank 1 $\Q$-line bundle $F$ on a reduced S2 projective variety $X$, they defined the $\mathbb{G}_m$-space
\[\mathcal{P}(F) \coloneqq \underline{\Spec}_X\bigoplus_{n\in\Z}F^{[n]}\] 
and also its $\mathbb{G}_m$-quotient stack \[\mathcal{X}_F\coloneqq [\mathcal{P}(F)/\mathbb{G}_m]\xrightarrow{\pi}X\]
One notes that there is a tautologically associated line bundle $L$ on it due to its map to $B\mathbb{G}_m$, and $\mathcal{X}_F$ is a stack uniformized by $L$. We call this construction the Abramovich--Hassett (AH) construction, and denote $\mathcal{X}_F$ by $AH(X,F)$. Hacking's canonical covering stack $\mathcal{X}$ agrees with $AH(X,\omega_X)$. In \cite{AH}, Abramovich and Hassett also proved that along the map $\pi$, for any $n\in\Z$,
\[\pi_*(L^{\otimes n}) = F^{[n]}\]
In fact, $AH(X,F)$ is the universal object that makes $F$ locally free (see Proposition \ref{universalprop}). In other words, the AH construction $\mathcal{X}_F\to X$ ``turns'' the $\Q$-Cartier divisor or divisorial sheaf into Cartier ones. Furthermore, assuming positivity on $F$ (a high enough reflexive power of $F$ being very ample, or equivalently, a high power of $L$ descends to a very ample line bundle on the coarse space), Abramovich and Hassett showed that one could normally embed the AH stack into a weighted projective stack $\bm{\mathcal{P}}(\rho)$ where $\rho$ is a $\mathbb{G}_m$-representation. But this is remarkable, because then they could construct the Hilbert scheme of all such AH stacks in $\bm{\mathcal{P}}(\rho)$, and then by some surgery one gets the moduli space $\mathcal{K}^{L}$ of so called ``uniformized twisted varieties'' (or taking rigidification $\mathcal{K}^{\lambda}=\mathcal{K}^L\rigidify \mathbb{G}_m$). This category is equivalent to the Koll\'ar category of $\Q$-line bundles:

\begin{defn*} [c.f. 5.2.1 of \cite{AH}, Koll\'ar family of $\Q$-line bundles]
    A Koll\'ar family of $\Q$-line bundles consists of the following information:\begin{itemize}
        \item $f:X\to B$ is a flat family of equi-dimensional (connected) reduced S2 schemes;
        \item $F\in\text{Coh}(X)$, such that\begin{itemize}
            \item [(1)] For each fiber $X_b$, the restriction $F|_{X_b}$ are reflexive of rank 1
            \item [(2)] For any $n\in\Z$, the formation of $F^{[n]}$ commutes with arbitrary base change $B'\to B$. This is know as the Koll\'ar's condition
            \item [(3)] For each geometric point $b$ of $B$, there exists an integer $N_{b} \neq 0$ such that $F|_{X_b}^{[N_{b}]}$ is invertible.
        \end{itemize}
    \end{itemize}
    \end{defn*}
\begin{thm*}
    [c.f. 5.3.6 of \cite{AH}] The category of Koll\'ar families of $\mathbb{Q}$-line bundles is equivalent to the category of uniformized twisted varieties via the base preserving functors
$$(X \rightarrow B, F)  \mapsto \left(\mathcal{X}_F \rightarrow B, L_F\right)$$
with $\mathcal{X}_F=\left[\mathcal{P}_F / \mathbb{G}_m\right]$, and its inverse
$$(\mathcal{X} \rightarrow B, L)  \mapsto \left(X \rightarrow B, \pi_* L\right)$$
where $X$ is the coarse moduli space of $\mathcal{X}$.
\end{thm*}
\noindent In particular, this equivalence immediately implies the isomorphism
$$\text{Def}_{\text{Kol}}(X,F) = \text{Def}(\mathcal{X},L)$$
where left hand side denotes the deformation space of $X$ together with a $\Q$-line bundle $F$ that satisfies Koll\'ar's conditions.
\subsection{About this paper} In much of the literature on the Minimal Model Program, KSBA stability, and birational geometry, the divisor $D$ in the pair $(X,D)$ is not just one marked divisor. $D$ is often a sum $D=\sum a_iD_i$ where the coefficients $a_i\in \Q$ (or in $\R$), and each $D_i$ serves as a $\Q$-Cartier marked divisor (see, for example, the discussion of KSBA stable families in \cite{kol23}). Therefore, a very natural question to ask is what if one wants to study the deformation theory where each marked divisor $D_i$ is allowed to deform separately. To be more precise, we ask the following questions: \begin{itemize}
    \item [(1)] Can we construct a ``multi version'' of the Abramovich--Hassett stack that turns all marked divisors $D_i$ into Cartier divisors? 
    \end{itemize}
    Assuming that we have a reduced S2 scheme $X$ and $F_1,F_2$ are two rank 1 $\Q$-Cartier reflexive sheaves, there are three very natural interpretations of the terminology ``doubly Abramovich--Hassett construction''. The first would be taking the relative Spec over $X$ of the direct sum of all reflexive combinations of $F_1$ and $F_2$, and then quotient by $\mathbb{G}_m^2$:
    \[\left[\underline{\Spec}_X\bigoplus_{(n_1,\,n_2)\in\Z^2}F_1^{[n_1]}[\otimes]_{\OO_X}F_2^{[n_2]}\middle/\mathbb{G}_m^2\right]\]
    The second is to first construct the usual AH stack using $F_1$: $\pi:AH(X,F_1)\to X$ and then construct on $AH(X,F_1)$ the AH stack using the reflexive pullback $(\pi^*F_2)^{[1]}$:
    \[AH(AH(X,F_1),(\pi^*F_2)^{[1]})\]
    or the other way round: use $F_2$ first and then the reflexive pullback of $F_1$. \`A priori they are not guaranteed to be the same, but in the following theorem, we prove that all these constructions agree. In particular, there is no ambiguity when we talk about the multi Abramovich--Hassett stacks:
\begin{thm} [= \ref{AHconstruction}]
Let $X$ be a reduced S2 scheme, $F_1,F_2$ be rank 1 reflexive $\OO_X$-modules that are both $\Q$-Cartier with Cartier indices $N_1,N_2$ respectively. Then the two possible doubly Abramovich--Hassett constructions are equivalent, i.e.
    $$AH(AH(X,F_1),(\pi^*F_2)^{[1]})\cong\left[\underline{\Spec}_X\bigoplus_{(n_1,n_2)\in\Z^2}F_1^{[n_1]}[\otimes]_{\OO_X}F_2^{[n_2]}\middle/\mathbb{G}_m^2\right]$$
    where $\pi:AH(X,F_1)\to X$ is the coarse space map. This holds more generally for $k$-fold AH construction where the left hand side is produced iteratively and the right hand side is \[\left[\underline{\Spec}_X\bigoplus_{(n_1,n_2,\,...,\,n_k)\in\Z^2}F_1^{[n_1]}[\otimes]_{\OO_X} \cdots [\otimes]_{\OO_X}F_k^{[n_k]}\middle/\mathbb{G}_m^k\right]\]
\end{thm}
 \begin{itemize}
     \item [(2)] Why do we want to look at the stack deformation space? What does it correspond to when passing to the coarse space?
 \end{itemize}   
In Section 3, we explain that the outcome of a $k$-fold AH construction is a ``$k$-uniformized twisted variety''; in particular, it is a globally $k$-diagonalizable (in short, $k$-cyclotomic) orbispace. The category of such stacks is exactly equivalent to the Koll\'ar category of $k$ $\Q$-line bundles:
\begin{thm}[= \ref{equivcat}]
    The category of Koll\'ar families of $n$ $\Q$-line bundles is equivalent to the category of $n$-uniformized twisted varieties via the base preserving functors 
    $$(X\to B,F_1,\,...,\,F_n)\mapsto (\mathcal{X}_{F_1,\,...,\,F_n}\to B,L_1,\,...,\,L_n)$$
    where $\mathcal{X}_{F_1,\,...,\,F_n}=[\mathcal{P}_{F_1,\,...,\,F_n}/\mathbb{G}_m^n]$ and $L_1,...,L_n$ are naturally associated line bundles on $\mathcal{X}$. Its inverse is simply the coarse moduli map
    $$(\mathcal{X}\to B,L_1,\,...,\,L_n)\mapsto (X\to B,\pi_*L_1,\,...,\,\pi_*L_n).$$
\end{thm}
\noindent This generalizes Theorem 5.3.6 of \cite{AH}. The stack deformation ensures that on the coarse space level, the $\Q$-line bundles still deform to $\Q$-line bundles, and they satisfy the natural generalization of Koll\'ar's condition:
\[\text{the formation of }F_1^{[a_1]}[\otimes]\cdots[\otimes]F_n^{[a_n]}\text{ for every $(a_i)_{i=1}^n\in \Z^n$ commutes with arbitrary base change}\]
\begin{itemize}
    \item [(3)] Can we embed multi-uniformized stack into a higher dimensional analogue of a weighted projective stack, and hence build the Hilbert scheme to construct moduli space?
\end{itemize}
The recent paper by Bragg, Olsson and Webb \cite{BOW} explained in their Section 4 how to embed a tame stack $\mathcal{X}$ with a det-ample vector bundle ($GL_k$-bundle) into a generalized Grassmannian stack $QP(V)$ where $V$ is a finitely generated module over the base ring with a $GL_k$-action. In Section 4 of this paper, we use their technology to embed our multi-uniformized stack into a toric stack (depending on a $\mathbb{G}_m^k$-representation). In our case, we have a $\mathbb{G}_m^k$-bundle, which is a direct sum of tautological line bundles raised to different powers. These powers are exactly given by the direction of polarization. 
\begin{thm} [= \ref{embedding}]
    Let $f:\mathcal{X}\to B$ be a proper $k$-cyclotomic stack, uniformized by $L_1,...,L_k$ (map $\mathcal{X}\to B\mathbb{G}_m^k$ is representable) with coarse space map $\pi:\mathcal{X}\to X$ and coarse space family $\bar{f}:X\to B$. If $L_1^{\otimes\alpha_1}\otimes L_2^{\otimes\alpha_2}\otimes\cdots\otimes L_k^{\otimes\alpha_k}$ descends to a very ample line bundle on $X$ (the $i$-th component of $\mathbb{G}_m^k$ acts on $L_i$ only), then we have a closed immersion $\mathcal{X}\hookrightarrow \QP(V)$ into a toric stack. Here $V$ is a direct sum of global section spaces of some $\Z^k$-combination of the $L_i$'s and depends on $\underline{\alpha}_i$'s.
\end{thm}
\noindent In Section 5, we do surgery on the Hilbert scheme to derive the algebraic stack $\operatorname{Sta}^{L_1,...,L_k,\,\underline{\alpha}}_{\mathfrak{F}}$, parametrizing $k$-uniformized $\underline{\alpha}$-normally embedded (see \ref{normalembedding}) stacks with Hilbert function $\mathfrak{F}$ and we show that
\begin{thm} [= \ref{stastack}]
    $\operatorname{Sta}^{L_1,\,...,\,L_k,\,\underline{\alpha}}_{\mathfrak{F}}$ is algebraic and locally of finite type.
\end{thm}
\noindent Further variants on this stack can give more specific moduli spaces of interest, e.g. the moduli of Seifert $\mathbb{G}_m^d$-bundles in Section 6.
\begin{itemize}
    \item [(4)] Can we further generalize the multi AH (mAH) construction?
\end{itemize}
One may have noticed that in order to perform the mAH construction, we don't necessarily need the components of the $\OO_X$-algebra to be reflexive combinations of $k$ rank 1 reflexive sheaves $F_1,...,F_k$. All we really need is a bunch of sheaves $\{\Ln\mid n\in\Z^d\}$ with a compatible multiplication structure on $\bigoplus_{\underline{n}\in\Z^d}\Ln$. Koll\'ar first studied the relative Spec of the direct sum of these sheaves under the name Seifert $\mathbb{G}_m$-bundles when the scheme is normal in \cite{Seifert}, and then discussed the semi-normal case in \cite{kol23}. Here the $\mathbb{G}_m^d$-action comes from the $\Z^d$-grading of the algebra. Therefore, we will call this category the Koll\'ar category of Seifert $\mathbb{G}_m^d$-bundles. The Koll\'ar's condition in this context would be:
\[\text{under arbitrary base change, the pullback of any }\Ln\text{ is reflexive on fibers.}\]
On the other hand, we remark that the key difference between the stacks constructed from these sheaves $\{\Ln\}$ and the usual multi AH stacks is that they are no longer $d$-uniformized twisted varieties (see 5.3.4 of \cite{AH}) because the coarse space map should \textit{not} necessarily be an isomorphism in codimension 1. The main theorem in Section 6 is:
\begin{thm} [= \ref{seiferteq}]
    The category of Koll\'ar families of Seifert $\mathbb{G}_m^d$-bundles is equivalent to the category of $d$-uniformized $d$-cyclotomic orbispace via the base preserving functors
    \[(X\to B,\{\Ln\}_{\underline{n}\in\Z^d})\mapsto (\mathcal{X}_{\{\Ln\}},\mathcal{L}_1,\,...,\,\mathcal{L}_d)\] where $\mathcal{X}_{\{\Ln\}}=[\mathcal{P}_{\{\Ln\}} /\mathbb{G}_m^d]$, and its inverse is given by the coase moduli map
    \[(\mathcal{X}\to B,\mathcal{L}_1,\,...,\, \mathcal{L}_d)\mapsto \left(X\to B, \Ln=(\pi_*\mathcal{L}_1^{\otimes n_1})[\otimes]\cdots[\otimes](\pi_*\mathcal{L}_d^{\otimes n_d})\right)\]
\end{thm}
\noindent and then we construct the moduli of Seifert $\mathbb{G}_m^d$-bundles:
\begin{prop} [= \ref{modulisft}]
    The moduli stack of all Seifert $\mathbb{G}_m^d$-bundles $(X\to B,\{\Ln\})$ with $\chi(X_b,\Ln)=\mathfrak{F}(\underline{n})$ for any $b\in B$ and for all $\underline{n}\in\Z^d$ is algebraic and locally of finite type.
\end{prop}
\noindent We end the paper by discussing the nice behavior when we stand on a family of normal varieties.

\subsection{Future directions} One other important motivation for this project was the moduli of AFI $\Q$-stable pairs $(X,D)$ that Filipazzi and Inchiostro introduced in \cite{AFI} (see also 8.24-26 of \cite{kol23}, here ``A'' in ``AFI'' stands for Alexeev). This notion of stability with floating coefficients guarantees that on the canonical model, divisors remains $\Q$-Cartier. In an upcoming paper with Bejleri \cite{bp}, we discuss how multi-uniformized stacks (or multi Abramovich--Hassett constructions) are used to understand the deformation of AFI stable families, and to construct the AFI moduli space
$$\mathcal{M}^{\operatorname{AFI}}_{n,\,p(t),\, I}$$
where $I\subset (0,1]\cap\Q$ provides coefficients for the boundary divisor $D=\sum a_iD_i$, $n=\dim X$ and $p(t) =(K_X+tD)^{n}$ is a fixed polynomial. The stability condition here is to impose that a multiple of $K_{X/S}+D$ is relatively nef and a multiple of $(K_{X/S}+(1-t)D)$ is relatively ample for any $t<\varepsilon_0$ small enough, rather than requiring a multiple of $K_{X/S}+D$ to be relatively ample.\medskip\\
\textbf{Conventions.} Throughout this paper, we work generally over arbitrary characteristic; and the families are over $\Spec \Z$. All the stacks involved will be tame stacks, which is in characteristic zero equivalent to the Deligne--Mumford (DM) condition (\cite{AOV}).
\medskip\\
\textbf{Acknowledgement.} The author expresses profound gratitude to his advisor, Dori Bejleri, for suggesting this project and for many insightful discussions. He also thanks his advisor at Harvard, Joe Harris, for continued guidance and support. The author is grateful to Dan Abramovich and J\'anos Koll\'ar for their helpful feedback on earlier drafts, and to Giovanni Inchiostro, Taeuk Nam and Anh Duc Vo for stimulating conversations. This project is partially funded by NSF grant DMS-2401483, and the author thanks the UMD Math Department and the Brin Center for their hospitality while part of this research was conducted.
\section{Construction of Doubly-Uniformized Stacks}
\noindent \textbf{Setup.} In this section, we work with a single scheme $X$ (over $\Spec\, k$) rather than with a family. $X$ is assumed to be a reduced scheme satisfying Serre's S2 condition. All sheaves in this paper will be Noetherian. A sheaf $F$ is called \textit{reflexive} if the natural map $F\to F^{\vee\vee}$ to its reflexive hull is an isomorphism. If $F$ is a rank 1 reflexive sheaf, then it's said to be \textit{$\Q$-Cartier} or \textit{$\Q$-line bundle} if $F^{[n]}$ is invertible for some $n\in\Z_{>0}$. The minimal such positive $n$ is called the \textit{Cartier index} of $F$. These notions still make sense when $X$ is a tame stack by imposing the conditions on the pullback of $F$ to an atlas. We use the following notations for singly AH construction:
\begin{align*}
    \mathcal{P}(F) &\coloneqq \underline{\Spec}_X \bigoplus_{n\in\Z} F^{[n]} \\
    AH(X,F) &= \mathcal{X}_F \coloneqq [\mathcal{P}(F)/\mathbb{G}_m]
\end{align*}
Due to the natural map $\mathcal{X}_F\to B\mathbb{G}_m$, we have an associated line bundle on $\mathcal{X}_F$, which we usually denote as $L$ or $L_F$. We call an open subset $U$ inside $X$ \textit{big} if the complement $X\setminus U$ has codimension $\geq2$. For clarity of argument, we prove all statements in the case $k=2$ (doubly AH), but everything naturally extends when $k\geq 3$.
\begin{prop}
    [Universal property of the Abramovich--Hassett construction] Let $X$ be a tame stack, $F$ be a reflexive rank 1 sheaf on $X$ that is $\Q$-Cartier. Let $\mathcal{Y}$ be a tame stack with a map $f:\mathcal{Y}\to X$ and a line bundle $\mathcal{G}$ on $\mathcal{Y}$ such that for every $n\in\Z$, $f_*\mathcal{G}^n\cong F^{[n]}$. Then $(\alpha:\mathcal{X}=\mathcal{X}_F\to X,L=L_F)$ is the universal object among $(f:\mathcal{Y}\to X,\mathcal{G})$, where $L$ comes from the $\mathbb{G}_m$-quotient on the AH construction $\mathcal{X}$.
\label{universalprop}
\end{prop}
\begin{pf}
    Since $\bigoplus_{n\in\Z} f_*{\mathcal{G}^n}=\bigoplus _{n\in\Z}F^{[n]}$ and by adjunction
    $$\hom \left(\bigoplus F^{[n]},\bigoplus f_*\mathcal{G}^n\right)=\hom\left(f^*\bigoplus F^{[n]},\bigoplus \mathcal{G}^n \right)$$
    we get a map $\mathcal{P}(\mathcal{G})\to \mathcal{P}(F)\times_X\mathcal{Y}\to\mathcal{P}(F)$. Now this map descends to
    $$[\mathcal{P}(\mathcal{G})/\mathbb{G}_m]\to [\mathcal{P}(F)/\mathbb{G}_m]$$
    because $[\mathcal{P}(\mathcal{G})/\mathbb{G}_m](T)$ contains the data of
    \[\begin{tikzcd}
	Z && {\mathcal{P}(\mathcal{G})} \\
	T
	\arrow["{\mathbb{G}_m\text{-equivariant}}", from=1-1, to=1-3]
	\arrow["{\mathbb{G}_m}"', from=1-1, to=2-1]
\end{tikzcd}\]
but $\mathcal{P}(\mathcal{G})\to\mathcal{P}(F)$ comes from a $\Z$-graded homomorphism, so it's also $\mathbb{G}_m$-equivariant. We get
\[\begin{tikzcd}
	Z && {\mathcal{P}(\mathcal{G})} & {\mathcal{P}(F)} \\
	T
	\arrow["{\mathbb{G}_m\text{-equivariant}}", from=1-1, to=1-3]
	\arrow["{\mathbb{G}_m}"', from=1-1, to=2-1]
	\arrow[from=1-3, to=1-4]
\end{tikzcd}\]
This defines an object in $[\mathcal{P}(F)/\mathbb{G}_m]$. This construction is clearly functorial. Now since $\mathcal{G}$ is a line bundle on $\mathcal{Y}$, $[\mathcal{P}(\mathcal{G})/\mathbb{G}_m]$ is simply $\mathcal{Y}$ and the right hand side is AH$(X,F)$.\medskip\\
Conversely, if we are given a map $\mathcal{Y}\to [\mathcal{P}(F)/\mathbb{G}_m]$, pulling back the associated line bundle on $[\mathcal{P}(F)/\mathbb{G}_m]$ gives a line bundle on $\mathcal{Y}$. This is the inverse of the previous construction. \qed
\end{pf}
\begin{thm} Let $X$ be a reduced S2 scheme, $F_1,F_2$ be rank 1 reflexive $\OO_X$-modules that are both $\Q$-Cartier with Cartier indices $N_1,N_2$ respectively. Then the two possible doubly Abramovich--Hassett constructions are equivalent, i.e.
    $$AH(AH(X,F_1),(\pi^*F_2)^{[1]})\cong\left[\underline{\Spec}_X\bigoplus_{(n_1,n_2)\in\Z^2}F_1^{[n_1]}[\otimes]_{\OO_X}F_2^{[n_2]}\middle/\mathbb{G}_m^2\right]=[\mathcal{P}({F_1,F_2})/\mathbb{G}_m^2]$$
    where $\pi:AH(X,F_1)\to X$ is the coarse moduli space map. We denote it by $AH(X,F_1,F_2)$ or $\mathcal{X}_{F_1,F_2}$.
    \label{AHconstruction}
\end{thm}
\begin{pf} We first locally take affine patch $X=\Spec A$ and let $F_1=B,F_2=C$ be reflexive rank 1 $A$-modules such that $B^{[N_1]}\cong A$ and $C^{[N_2]}\cong A$. We have the following maps and affine patches:
\[\begin{tikzcd}
	{\bigoplus_{n=0}^{N_1-1}B^{[n]}} & {\underline{\text{Spec}}_X\bigoplus_{n=0}^{N_1-1}F_1^{[n]}} \\
	{\bigoplus_{n=0}^{N_1-1}B^{[n]}+\mu_{N_1}\text{-action}} & {X_1=\left[\underline{\text{Spec}}_X\bigoplus_{n=0}^{N_1-1}F_1^{[n]}/\mu_{N_1}\right]=\left[\underline{\text{Spec}}_X\bigoplus_{n\in\Z}F_1^{[n]}/\mathbb{G}_m\right]} \\
	A & X
	\arrow["v", from=1-2, to=2-2]
	\arrow["\pi", from=2-2, to=3-2]
\end{tikzcd}\]
where the two expressions for $X_1$ are the same due to Corollary 4.4 of \cite{BejInchio}.\medskip\\
By descent, we have that for each $m\in\Z$,
\begin{align*}(\pi^*F_2)^{[m]} &=\left(\bigoplus_{n=0}^{N_1-1}B^{[n]}\otimes_A C^{\otimes m}\right)^{\vee\vee} \text{ + grading on $n$ (reflexive hull taken w.r.t. $\oplus_nB^{[n]}$)}
\end{align*}
Write $R=\bigoplus_{n=0}^{N_1-1} B^{[n]}$ an $A$-algebra, and $T=C^{\otimes m}$ an $A$-module. $\Spec A$ is an S2 scheme; $T$ is a rank 1 $A$-module, and $R$ is a reflexive $A$-module since finite direct sum commutes with taking dual. It's worth noting that here we must take the \'etale cover $\underline{\text{Spec}}_X\bigoplus_{n=0}^{N_1-1}F_1^{[n]}$ instead of the smooth cover $\underline{\text{Spec}}_X\bigoplus_{n\in\Z}F_1^{[n]}$ because then $R=\bigoplus_{n\in\Z}B^{[n]}$ is not reflexive: the dual of infinite direct sum becomes infinite product.
\begin{lem} [A]
    $\hom_A(T,R)\cong \hom_A(T,A)[\otimes]_AR$ as $A$-modules, and as $R$-modules.
\end{lem}
\begin{pf}
    Consider the obvious map $\hom_A(T,A)\otimes_AR\xrightarrow{\varphi} \hom_A(T,R)$. Away from a codimension 2 locus $Z\subset\Spec A$, $T$ and $R$ are both locally free, so $\varphi$ is an isomorphism away from $Z$. Since $R$ is reflexive, we have
    \begin{align*}
        \hom_A(T,R)&=\hom_A(T,\hom_A(R^{\vee},A))\\
        &= \hom_A(T\otimes_A R^{\vee},A)\\
        &= (T\otimes_A R^{\vee})^{\vee}
    \end{align*}
    so $\hom_A(T,R)$ is reflexive. Now $\varphi$ is an isomorphism away from codimension 2, $A$ is S2, and the target of $\varphi$ is reflexive, so by the universal property we get isomorphism 
    $$(\hom_A(T,A)\otimes_AR)^{\vee\vee}\xrightarrow{\cong} \hom_A(T,R)$$
    where the reflexive hull is taken with respect to $A$ and this is an $R$-module isomorphism as well. \qed
\end{pf}
\begin{lem}
    [B] $(\pi^*F_2)^{[m]}\cong \bigoplus_{n=0}^{N_1-1}B^{[n]}[\otimes]_A C^{[m]}.$
\end{lem}
\begin{pf}
        Let $U\subset \Spec A$ be the big open subset where $B,C$ are locally free. Let $V=(p=v\circ\pi)^{-1}(U)$. Let $i:U\hookrightarrow \Spec A$ and $j:V\hookrightarrow \Spec R$. Call $M=\hom_A(T,A)$.
        \[\begin{tikzcd}
	V & {\Spec R} \\
	U & {\Spec A}
	\arrow["j", hook, from=1-1, to=1-2]
	\arrow["q"', from=1-1, to=2-1]
	\arrow["p", from=1-2, to=2-2]
	\arrow["i", hook, from=2-1, to=2-2]
\end{tikzcd}\]
As an $A$-module, 
\begin{align*}
    R[\otimes]_AM &= i_*i^* p_* p^* M\\
    &= i_*q_*j^*p^*M \qquad\qquad\text{[by push-pull]}\\
    &= (p\circ j)_*(p\circ j)^*M
\end{align*}
thus, as an $R$-module, $R[\otimes]_AM$ is simply $j_*(p\circ j)^*M$. Now 
\begin{align*}
    \hom_R(\hom_R(R\otimes_A T,R),R)&=\hom_R(\hom_A(T,R),R)\qquad\qquad
     &\text{[adjunction]}\\
    &= \hom_R(R[\otimes]_A\hom_A(T,A),R) &\text{[by lemma A]}\\
    &= \hom_R(j_*(p\circ j)^*M,R)\\
    &= \hom_R(j_*j^*p^*M,R)\\
    &= \hom_V((p\circ j)^*M,\OO_V) &[*]\\
    &= \hom_A(M,p_*j_*\OO_V) &\text{[adjunction]}\\
    &= \hom_A(M,R)\\
    &= R[\otimes]_A\hom_A(M,A) &\text{[by lemma A]}\\
    &= R [\otimes]_A \hom_A(\hom_A(T,A),A) 
    \end{align*}
    Here $[*]$ is due to the fact that $\mathcal{H}\text{om}_R(j_*j^*p^*M,R)$ is reflexive with respect to $R$, so it equals $j_*\mathcal{H}\text{om}_V(j^*p^*M,\OO_V)$. This also relies on the fact that $R$ is S2 and $V$ is a big open in $\Spec R$. Note that in this proof we use finitely generated modules and coherent sheaves interchangeably.\qed
\end{pf}
\noindent The next thing to do is to construct relative Spec on the AH stack $X_1$. Since we had an \'etale cover $\bigoplus_{n=0}^{N_1-1}B^{[n]}$ on $X_1$, by descent, we perform the construction on $\bigoplus_{n=0}^{N_1-1}B^{[n]}$ and take the $\mu_{N_1}$ quotient induced by the grading on $n$:
$$\bigoplus_{m\in\Z}\bigoplus_{n=0}^{N_1-1}B^{[n]}[\otimes]_A C^{[m]}\quad\text{with $\mu_{N_1}$ action induced by grading on $n$}$$
Thus, by further quotienting the $\mathbb{G}_m$-action induced by grading on $m$, 
\begin{align*}
    AH(AH(X,F_1),(\pi^*F_2)^{[1]}) &= \left[\underline{\Spec}_X\bigoplus_{m\in\Z}\bigoplus_{n=0}^{N_1-1
}B^{[n]}[\otimes]_A C^{[m]}/\mu_{N_1} \middle/\mathbb{G}_m\right]\\
    &=\left[\underline{\Spec}_X\bigoplus_{m=0}^{N_2-1}\bigoplus_{n=0}^{N_1-1
}B^{[n]}[\otimes]_A C^{[m]}\middle/\mu_{N_1}\times\mu_{N_2}\right]
\end{align*}
\begin{rem}
    Under our construction, if $\bigoplus_{n\in\Z} B^{[n]}[\otimes]_A C^{[N']}$ happens to be locally free, i.e. $\cong \bigoplus_{n\in\Z}B^{[n]}$, for some $N'|N_2$, then we could simply take a further quotient by $\mu_{N_2}/\mu_{N'}$, to obtain
    $$\left[\underline{\Spec}_X\bigoplus_{m=0}^{N'-1}\bigoplus_{n=0}^{N_1-1
}B^{[n]}[\otimes]_A C^{[m]}\middle/\mu_{N_1}\times\mu_{N'}\right]$$
In particular, when $N'=1$, the second AH construction is redundant.
\end{rem}
\begin{rem}
    If we let $F_2=F_1$, then since
    $$\bigoplus_{n\in\Z}B^{[n]}[\otimes]_A B^{[1]}\cong_{\oplus B^{[n]}\text{-mod}}\bigoplus_{n\in\Z}B^{[n]}$$
which means that $(\pi^*F_2)^{[1]}$ is locally free, we don't need the second AH construction.
\end{rem}
\noindent Now we construct a global map 
$$\left[\underline{\Spec}_X\bigoplus_{\Z^2}F_1^{[n_1]}[\otimes]_{\OO_X} F_2^{[n_2]}\middle/\mathbb{G}_m^2\right]\to AH\left(AH(X,F_1),(\pi^*F_2)^{[1]}\right)$$
so that the previous local computation proves the claimed isomorphism.
\[\begin{tikzcd}
	& {X_2=\left[\underline{\Spec}_{X_1}\bigoplus_{m\in\Z}(\pi^*F_2)^{[m]}/\mathbb{G}_m\right]} \\
	{Y=\left[\underline{\Spec}_X\bigoplus_{\Z^2}F_1^{[n_1]}[\otimes]_{\OO_X} F_2^{[n_2]}/\mathbb{G}_m^2\right]} & {X_1=\left[\underline{\Spec}_X\bigoplus_{n\in\Z}F_1^{[n]}/\mathbb{G}_m\right]} \\
	& X
	\arrow["{\pi'}", from=1-2, to=2-2]
	\arrow["\gamma", dashed, from=2-1, to=1-2]
	\arrow["\beta"', dashed, from=2-1, to=2-2]
	\arrow["\alpha"', from=2-1, to=3-2]
	\arrow["\pi", from=2-2, to=3-2]
\end{tikzcd}\]
Let $\alpha,\pi,\pi'$ be the obvious maps as shown in the diagram. Here $Y$ as a $\mathbb{G}_m^2$-quotient tautologically carries two line bundles $L_1$ and $L_2$. $(\alpha^*F_1)^{[1]}$ is equal to the Cartier divisor $L_1$ on $Y$, so $\alpha_*L_1^n=F_1^{[n]}$ for any $n\in\Z$. Hence we get the map $\beta$ due to universal property. Since $\beta^*(\pi^*F_2)^{[1]}$ is equal to the Cartier divisor $L_2$ on $Y$, so $\beta_*L_2^n=(\pi^*F_2)^{[n]}$ for any $n\in\Z$, and we get the map $\gamma$ by the universal property. This concludes the proof of the theorem.\qed\end{pf}
\begin{prop}
    [Universal property of doubly AH stacks] Let $X$ be a scheme and $F_1,F_2$ be two reflexive rank 1 sheaves that are $\Q$-Cartier. Let $\mathcal{X}=AH(X,F_1,F_2)$ be the doubly AH stack, uniformized by line bundles $L_1$ and $L_2$. Then $(\mathcal{X}\to X,L_1,L_2)$ is the universal object among $(f:\mathcal{Y}\to X,\mathcal{G}_1,\mathcal{G}_2)$ where $\mathcal{Y}$ is a tame stack, $\mathcal{G}_1,\mathcal{G}_2$ are two line bundles on $\mathcal{Y}$ such that $f_*(\mathcal{G}_1^n\otimes \mathcal{G}_2^m)\cong F_1^{[n]}[\otimes]F_2^{[m]}$.
\end{prop}
\begin{pf}
    This is analogous to the singly AH case. We use the isomorphism $\bigoplus f_*(\mathcal{G}_1^n\otimes\mathcal{G}_2^m)\cong\bigoplus F_1^{[n]}[\otimes]F_2^{[m]}$ and adjunction $$\hom\left(\bigoplus F_1^{[n]}[\otimes]F_2^{[m]}, f_*(\mathcal{G}_1^n\otimes\mathcal{G}_2^m)\right)=\hom \left(f^*\bigoplus F_1^{[n]}[\otimes]F_2^{[m]},\mathcal{G}_1^n\otimes\mathcal{G}_2^m\right)$$
    to get $\mathcal{P}(\mathcal{G}_1,\mathcal{G}_2)\to \mathcal{P}(F_1,F_2)\times_X\mathcal{Y}\to\mathcal{P}(F_1,F_2)$. This descends to \[[\mathcal{P}(\mathcal{G}_1,\mathcal{G}_2)/\mathbb{G}_m^2]\to[\mathcal{P}(F_1,F_2)/\mathbb{G}_m^2]\] But this is just $\mathcal{Y}\to \mathcal{X}$. Such map is unique due to the chosen line bundles $\mathcal{G}_1,\mathcal{G}_2$.\qed
\end{pf}
\begin{rem}
    The construction we had for doubly AH stacks works in general for $k$-fold AH: we have the following isomorphism:
    \[\hspace{-0.5cm}\left[\underline{\Spec}_X\bigoplus_{n_1,\,...,\,n_k\in\Z^k} F_1^{[n_1]}[\otimes] F_2^{[n_2]} [\otimes] \cdots [\otimes] F_k^{[n_k]}\middle/\mathbb{G}_m^k\right] =AH(AH(\cdots (AH(X,F_1),F_2^{[1]}),F_3^{[1]}),\cdots ),F_k^{[1]})\]
    which we denote by $\mathcal{X}$. Together with the line bundles $L_1,...,L_k$ given by the map to $B\mathbb{G}_m^k$, this is the universal object that makes every reflexive combination of $F_1,...,F_k$ invertible, in the sense of the previous proposition.
\end{rem}

\section{Koll\'ar Category and Doubly-Uniformized Twisted Varieties}
\begin{defn} [Koll\'ar families] \label{kolfam}
    By a Koll\'ar family of double $\Q$-line bundles, we mean:\begin{itemize}
        \item $f:X\to B$ is a flat family of equi-dimensional (connected) reduced S2 schemes;
        \item $F_1,F_2\in\text{Coh}(X)$, such that\begin{itemize}
            \item [(1)] For each fiber $X_b$, the restriction $F_1|_{X_b}, F_2|_{X_b}$ are reflexive of rank 1
            \item [(2)] For any $n,m\in\Z$, the formation of $F_1^{[n]}[\otimes]_{\OO_X}F_2^{[m]}$ commutes with arbitrary base change $B'\to B$
            \item [(3)] For each geometric point $s$ of $B$, there exist integers $N_{s,1},N_{s,2} \neq 0$ such that $F_1|_{X_s}^{[N_{s,1}]}$ and $F_2|_{X_s}^{[N_{s,2}]}$ are both invertible.
        \end{itemize}
    \end{itemize}
    A morphism from a Koll\'ar family $(X\to B,F_1,F_2)$ of double $\Q$-line bundles to another $(X'\to B',F_1',F_2')$ consists of a Cartesian diagram
    \[\begin{tikzcd}
	X & {X'} \\
	B & {B'}
	\arrow["\varphi", from=1-1, to=1-2]
	\arrow["f", from=1-1, to=2-1]
	\arrow["{f'}", from=1-2, to=2-2]
	\arrow[from=2-1, to=2-2]
\end{tikzcd}\] with isomorphisms $\alpha:F_1\to \varphi^*F_1'$ and $\beta:F_2\to\varphi^*F_2'$. These define the objects and morphisms in the category of Koll\'ar families of double $\Q$-line bundles.
\label{kolfam}
\end{defn}
\begin{rem}
    Condition (2) implies that $F_1^{[n]}[\otimes]F_2^{[m]}$ are flat over $B$ for all $m,n\in\Z$ due to theorem 10 in \cite{Kol95} (or proposition 5.1.4 in \cite{AH}). In particular, $F_1,F_2$ are flat over $B$. 
\end{rem}
\begin{lem}
    Let $(f:X\to B,F_1,F_2)$ be a Koll\'ar family of double $\Q$-line bundles. For $s\in B$ a geometric point, there exists $N_{s,i}\in \Z$ such that $F_i^{[N_{s,i}]}|_{X_s}$ is invertible (for $i=1,2$). Then there's a neighbourhood $U$ of $s$ in $B$ such that $F_1^{[N_{s,1}]}|_{X_U}$ and $F_2^{[N_{s,2}]}|_{X_U}$ are both invertible.
\end{lem}

\begin{pf}
    Analogous to the argument for Lemma 5.2.2 in \cite{AH}.\qed
\end{pf}
\begin{rem}
    Under condition (2) and Noetherian condition on the base $B$, $F_1$ or $F_2$ is $\Q$-Cartier if and only if it's $\Q$-Cartier when restricted to the geometric fibers.
\end{rem}
\begin{defn}
    The doubly AH stack for the Koll\'ar family $(X\to B,F_1,F_2)$ is
    $$\mathcal{X}_{F_1,F_2}\coloneqq\left[ \underline{\Spec}_X \bigoplus _{n,m\in\Z^2} F_1^{[n]} [\otimes]_{\OO_X} F_2^{[m]}\middle/\mathbb{G}_m^2\right]$$
    where $\mathbb{G}_m^2$-action is induced by gradings on $n,m$, and $L_1 ,L_2$ are associated line bundles. \label{dahforkol}
\end{defn}

\begin{prop}
    For a Koll\'ar family of $k$ $\Q$-line bundles $(X\to B,F_1,F_2)$ over a Noetherian base $B$, we still have the isomorphism 
    \[\mathcal{X}_{F_1,F_2}\cong AH(AH(X,F_1),(\pi^*F_2)^{[1]})\]
    where $\pi:AH(X,F_1)\to X$ is the AH stack formed by $F_1$. \label{kolfamAH}
\end{prop}
\begin{pf}
    Pick $N_1,N_2\in\Z_{>0}$ Cartier indices of $F_1,F_2$ that work for every fiber $X_b$. The proof of universal property still works in this case, so we get a map
    \[\gamma:\mathcal{X}_{F_1,F_2} \to AH(AH(X,F_1),(\pi^*F_2)^{[1]})\]
    and by the Koll\'ar's condition \ref{kolfam}(2), fiberwise AH constructions of both sides are the same as the AH constructions restricted to the fibers. Hence we can apply Theorem \ref{AHconstruction} on the fibers to see that $\gamma|_b$ are all isomorphisms, and hence $\gamma$ is an isomorphism.\qed
\end{pf}

\begin{rem}
    Definition \ref{kolfam} can be generalized to include $n$ $\Q$-line bundles to define the Koll\'ar family of $n$ $\Q$-line bundles $(f:X\to B,F_1,...,F_n)$. Condition (2) should be replaced by: for any $a_1,...,a_n\in\Z$, the formation of $$F_1^{[a_1]}[\otimes]\cdots [\otimes] F_n^{[a_n]}$$ commutes with arbitrary base change $B'\to B$. Same argument in the proof of Proposition \ref{kolfamAH} also implies that there's no ambiguity to talk about multi AH constructions over Koll\'ar families of $n$ $\Q$-line bundles.
\end{rem}

\begin{defn} [c.f. 4.1.1 of \cite{AH}]
    An orbispace is a stack $\mathcal{X}$ that is separated, with finite diagonal, of finite type over a field, equidimensional, geometrically connected and reduced, and contains an open dense $U\subset X$ that is an algebraic space.
\end{defn}
\begin{defn}
    A stack $\mathcal{X}$ is $n$-diagonalizable if there is a rank $n$ torus $\mathcal{T}\to\mathcal{X}$ such that the inertia stack $\mathcal{I}_\mathcal{X}$ embeds into $\mathcal{T}$ over $\mathcal{X}$. For any geometric point $ x\in\mathcal{X}$, the automorphism group $\text{Aut}_{\mathcal{X}}(x)$ is a finite subgroup of an $n$-torus.
    \[\begin{tikzcd}
	{\mathcal{I}_{\mathcal{X}}} && {\mathcal{T}} \\
	& {\mathcal{X}}
	\arrow[hook, from=1-1, to=1-3]
	\arrow[from=1-1, to=2-2]
	\arrow[from=1-3, to=2-2]
\end{tikzcd}\]
\end{defn}
\begin{defn}
    An $n$-cyclotomic stack or a globally $n$-diagonalizable stack (resp. orbispace) is a tame stack (resp. orbispace) $\mathcal{X}$ such that we can embed the inertia stack $\mathcal{I}_{\mathcal{X}}$ into a global rank $n$ torus $T\times \mathcal{X}$. For any geometric point $ x\in\mathcal{X}$, the automorphism group $\text{Aut}_{\mathcal{X}}(x)$ is a finite subgroup of the $n$-torus $T$.
    \[\begin{tikzcd}
	{\mathcal{I}_{\mathcal{X}}} && {T\times\mathcal{X}} \\
	& {\mathcal{X}}
	\arrow[hook, from=1-1, to=1-3]
	\arrow[from=1-1, to=2-2]
	\arrow[from=1-3, to=2-2]
\end{tikzcd}\] \label{ncyclotomic}
    \end{defn}
\begin{prop}
    [Functor] Let $\mathcal{X}_{F_1,F_2}\to X\to B$ be as in Definition \ref{dahforkol}, then
    \begin{itemize} 
        \item [(1)] $\mathcal{X}_{F_1,F_2}\to B$ is a family of bicyclotomic orbispaces with S2 fibers.
        \item [(2)] $\pi:\mathcal{X}_{F_1,F_2}\to X$ is the coarse moduli space map. It's an isomorphism on the open subset where $F$ is invertible, complement of which has codimension $>1$ in each fiber.
        \item [(3)] For any $e,f\in\Z$, $\pi_*\left(L_1^{e}\otimes L_2^{f}\right)=F_1^{[e]}[\otimes] F_2^{[f]}$.
        \item [(4)] Given a morphism in the Koll\'ar category $(\varphi,\alpha,\beta)$ from $(X\to B,F_1,F_2)$ to $(X'\to B',F_1',F_2')$, we canonically have isomorphisms $\mathcal{P}_{F_1,F_2}\cong\varphi^* \mathcal{P}_{F_1',F_2'}$ and $\mathcal{X}_{F_1,F_2}\cong\varphi^*\mathcal{X}_{F_1',F_2'}$.\label{functor}
    \end{itemize}
    Thus, we get a functor from the Koll\'ar category to the category of doubly-uniformized twisted varieties.
\end{prop}
\begin{pf}
    (1) working locally on $B$, due to the lemma above we may assume $F_1^{[N_1]}$ and $F_2^{[N_2]}$ are locally free for some $N_1,N_2>0$. Consider the natural homomorphisms
    $$m_{e,f}:\left(F_1^{[N_1]}\right)^{e}\otimes \left(F_2^{[N_2]}\right)^{f}\to F_1^{[N_1e]}[\otimes] F_2^{[N_2f]}$$
    Using S2 and the fact that $F_i$'s are locally free on a big open in each fiber $X_b$, we conclude that $m_{e,f}$ is an isomorphism. Consider the $\mathbb{G}_m^2$-equivariant map
    $$\mathcal{P}_{F_1,F_2}\to\underline{\Spec}_X\left(\bigoplus_{m,n\in\Z}F_1^{[nN_1]}[\otimes]F_2^{[mN_2]}\right)$$
    with $\mathbb{G}_m^2$ acts on the target with stabilizer being a subgroup of $\mu_{N_1}\times\mu_{N_2}$. Since this morphism is representable, the stabilizers on $\mathcal{P}_{F_1,F_2}$ are finite subgroups of the torus $\mathbb{G}_m^2$. A fiber in $\mathcal{P}_{F_1,F_2}$ over $b\in B$ is S2 as the spectrum of an algebra with reflexive components over an S2 base (here we are using the assumption that the formation of $F_1^{[n]}[\otimes]F_2^{[m]}$ commutes with arbitrary base change; in this case, base change to a point). $\mathcal{P}_{F_1,F_2}\to \mathcal{X}$ is smooth and surjective, so the quotient stack $\mathcal{X}_b$ is also S2. \medskip\\
    For (2), since the invariant part of $\bigoplus _{n,m\in\Z^2} F_1^{[n]} [\otimes]_{\OO_X} F_2^{[m]}$ is simply the degree $(0,0)$ piece, i.e. $\OO_X$, we conclude that $X$ is the coarse moduli space of $\mathcal{X}$. Over the locus where $F_1,F_2$ are both invertible, the scheme $\mathcal{P}_{F_1,F_2}$ is a principal $\mathbb{G}_m^2$-bundle, and so $\mathcal{X}|_U=\left[\mathcal{P}_{F_1,F_2}|_U\middle/\mathbb{G}_m^2\right]=U$.\medskip\\
    For (3), note that $L_1^e\otimes L_2^f$ is the degree $(e,f)$ component of the $\OO_{\mathcal{X}}$-algebra $\bigoplus L_1^{n}\otimes_{\OO_{\mathcal{X}}} L_2^{m}$, whose spectrum is $\mathcal{P}_{F_1,F_2}$, and the degree $(e,f)$ component of $\bigoplus F_1^{[n]} [\otimes]_{\OO_X} F_2^{[m]}$ is $F_1^{[e]} [\otimes]_{\OO_X} F_2^{[f]}$.\medskip\\
    For (4), again using the assumption that the formation of $F_1^{[n]}[\otimes]F_2^{[m]}$ commutes with arbitrary base change, $F_1^{[n]}[\otimes]_{\OO_X}F_2^{[m]}=\varphi^*\left( F_1'^{[n]}[\otimes]F_2'^{[m]}\right)$. Together with functoriality of taking relative Spec, we complete the proof.\qed
\end{pf}
\begin{defn}
    A family of doubly-uniformized twisted variety $(\mathcal{X}\to B,L_1,L_2)$ is a flat family of bicyclotomic orbispaces uniformized by $L_1,L_2$ (that means, the automorphism group at each point acts faithfully and diagonally on $L_1\oplus L_2$), such that $\pi:\mathcal{X}\to X$ is an isomorphism away from a subset of codimension $>1$ in each fiber.
\end{defn}
\begin{defn}
    More generally, a family of $n$-uniformized twisted variety $(\mathcal{X}\to B,L_1,...,L_n)$ is a flat family of $n$-cyclotomic orbispaces uniformized by $L_1,...,L_n$ (again, the automorphism group at each point acts faithfully and diagonally on $L_1\oplus\cdots\oplus L_n$), such that $\pi:\mathcal{X}\to X$ is an isomorphism away from a subset of codimension $>1$ in each fiber.
\end{defn}
\begin{rem}
   Notice here that the torus $T$ in Definition \ref{ncyclotomic} is simply $\text{Aut}(L_1)\times\cdots\times\text{Aut}(L_n).$
\end{rem}
\begin{prop}
    $n$ line bundles $L_1,...,L_n$ on $\mathcal{X}$ form uniformizing line bundles for $\mathcal{X}\to B$ if and only if the map $\mathcal{X}\to B\mathbb{G}_m^n$ induced by $L_1,...,L_n$ is representable.
\end{prop}
\begin{prop}
    If $v_1,...,v_n\in\Z^n$ is a $\Z$ basis for $\Z^n$, then an $n$-uniformized stack $(\mathcal{X},L_1,...,L_n)$ is also uniformized by $L^{v_1},...,L^{v_n}$, where $L^{v_i}$ means $L_1^{v_{i,1}}\otimes\cdots\otimes L_n^{v_{i,n}}$.
\end{prop}
\begin{rem}
    Note that we are still working over orbispaces even though leaving the cyclotomic world. So, as in the singly AH construction, 2-isomorphism is unique when it exists, and hence it's equivalent to its associated category, and morphisms are isomorphism classes of 1-morphisms.
\end{rem}
\noindent\textbf{Some lemmas on $n$-cyclotomic stacks.}
\begin{lem}
    $n$-cyclotomic stacks are tame.
\end{lem}
\begin{pf}
    Any $\mu_m$ is linearly reductive, and a subgroup of the product of $n$ linearly reductive groups is still linearly reductive. Hence by Theorem 3.2 (b) of \cite{AOV}, we conclude that any $n$-cyclotomic stack is tame.\qed
\end{pf}
\begin{lem}
    Consider an $n$-cyclotomic stack $f :\mathcal{X} \to B$ and a geometric closed point $\xi: \Spec K \to\mathcal{X}$ with stabilizer $G\subset \mathbb{G}_m^n$. $G$ is finite and is a product of at most $n$ $\mu_{m_i}$-groups. Let $\bar{\xi} : \Spec K \to X$ be the corresponding point. Then, in a suitable \'etale neighborhood of $\bar{\xi}$, the stack $\mathcal{X}$ is a quotient of an affine scheme by an action of $G$.
\end{lem}
\begin{pf}
    Follows from Theorem 3.2 (d) of \cite{AOV}.\qed
\end{pf}
\begin{lem}
    Let $f:\mathcal{X}\to B$ be a flat family of $n$-cyclotomic stacks. Then \begin{itemize}
        \item [(1)] The coarse moduli space $\bar{f}:X\to B$ is flat;
        \item [(2)] Formation of the coarse moduli space commutes with base extension, namely, for any morphism of schemes $X'\to X$, the coarse moduli space of $\mathcal{X}\times_X X'$ is $X'$. In particular, the geometric fibers of $\bar{f}$ are coarse moduli spaces of the corresponding fibers of f.
    \end{itemize}
\end{lem}
\begin{pf}
    Follows from Corollary 3.3 of \cite{AOV}.\qed
\end{pf}

\begin{thm}
    The category of Koll\'ar families of double $\Q$-line bundles is equivalent to the category of doubly-uniformized twisted varieties via the base preserving functors 
    $$(X\to B,F_1,F_2)\mapsto (\mathcal{X}_{F_1,F_2}\to B,L_1,L_2)$$
    where $\mathcal{X}_{F_1,F_2}=[\mathcal{P}_{F_1,F_2}/\mathbb{G}_m^2]$ and $L_1,L_2$ are naturally associated line bundles on $\mathcal{X}$. Its inverse is simply the coarse moduli map
    $$(\mathcal{X}\to B,L_1,L_2)\mapsto (X\to B,\pi_*L_1,\pi_*L_2).$$
\end{thm}
\begin{pf}
    The Proposition \ref{functor} above shows that there is a functor from the Koll\'ar category to the category of doubly-uniformized twisted varieties. Conversely, we take a doubly-uniformized twisted variety $(\mathcal{X}\to B,L_1,L_2)$, with coarse moduli $\pi:\mathcal{X}\to X$. The formation of coarse moduli space commutes with arbitrary base change, so the universal property of coarse moduli space guarantees that each fiber of $X\to B$ is reduced. If $N_1$ is the index of $\mathcal{I}_{\mathcal{X}}\hookrightarrow T\times\mathcal{X}$ with respect to the first component of the torus, $N_2$ is the index for the second component, then $L_i^{N_i}$ descends to an invertible sheaf on $X$ which coincides with $\pi_*L_i^{N_i}=(\pi_*L_i)^{[N_i]}$.\medskip\\ Also, any geometric point of $\mathcal{X}$ admits an \'etale neighbourhood isomorphic to the quotient of an S2 affine scheme $V $ over a linearly reductive group $G$, which is a product of $\mu_{m_i}$-groups. The local coarse moduli space is the GIT quotient $V/G$, which is still S2 since the invariants are direct summands in the ring $\OO_V$. Similarly, $\pi_*(L_1^j\otimes L_2^k)$ are direct summands in the algebra of $\mathcal{P}$, which is affine over $X$, flat over $B$ with S2 fibers. Thus these sheaves are flat over $B$, saturated, and their formation commutes with arbitrary base change $B'\to B$. \qed
\end{pf}
\noindent In general, by replacing 2 with $n$, we would get the equivalence between the category of Koll\'ar families of $n$ $\Q$-line bundles and the category of $n$-uniformized twisted varieties due to the same arguement:
\begin{thm}
    The category of Koll\'ar families of $n$ $\Q$-line bundles is equivalent to the category of $n$-uniformized twisted varieties via the base preserving functors 
    $$(X\to B,F_1,\,...,\,F_n)\mapsto (\mathcal{X}_{F_1,\,...,\,F_n}\to B,L_1,\,...,\,L_n)$$
    where $\mathcal{X}_{F_1,\,...,\,F_n}=[\mathcal{P}_{F_1,\,...,\,F_n}/\mathbb{G}_m^n]$ and $L_1,...,L_n$ are naturally associated line bundles on $\mathcal{X}$. Its inverse is simply the coarse moduli map
    $$(\mathcal{X}\to B,L_1,\,...,\,L_n)\mapsto (X\to B,\pi_*L_1,\,...,\,\pi_*L_n).$$
    \label{equivcat}
\end{thm}
\section{Embedding Doubly-Uniformized Bicyclotomic Stack into a Toric Stack}
\noindent Our goal in this section is to provide an embedding of a $k$-uniformized $k$-cyclotomic stack into a quasi-projective toric stack. With such embedding we could construct the relevant moduli spaces by first taking the corresponding Hilbert scheme. This is a higher rank analogue of Corollary 2.4.4 in \cite{AH} where Abramovich and Hassett embedded the singly AH stack into a weighted projective stack. For clarity of the arguments we continue to assume that $k=2$.\medskip\\
Let $f:\mathcal{X}\to B$ be a proper bicyclotomic stack that is uniformized by $L_1,L_2$ (i.e. the morphism $\mathcal{X}\to B\mathbb{G}_m^2$ given by $L_1,L_2$ is representable). We denote by $\pi:\mathcal{X}\to X$ the coarse moduli space map, and $\bar{f}:X\to B$ the induced coarse space family. 
\begin{defn}
    We say that $(\mathcal{X},L_1,L_2)$ is $(a,b)$-polarized if there is an $\bar{f}$-ample line bundle $M$ on $X$ and a positive integer $N$ such that $\pi^*M \cong (L_1^a\otimes L_2^b)^N$.
\end{defn}
\noindent Let $M_1,M_2$ be multiples of $L_1,L_2$ respectively that descend to the coarse space $X$ and we view them as line bundles on $X$. Assume that the relative ample cone for $\bar{f}$ meets the plane $\Q\langle M_1,M_2\rangle$. By openness of ample cone, there is a wedge $\Lambda$ in $\Q\langle L_1,L_2\rangle$ such that any point in the interior of $\Lambda$ is polarizing. It is always possible to find a $\Z$-basis $\{L_1',L_2'\}\subset\Lambda$ for $\Z\langle L_1,L_2\rangle$. Without loss of generality, we may assume $(\mathcal{X},L_1,L_2)$ is such that both $L_1$ and $L_2$ are already polarizing, i.e. both $M_1,M_2$ are $\bar{f}$-ample.\medskip\\\noindent
\textbf{Notation.} For the most generality, we will only make the assumption that $L_1^{\alpha}\otimes L_2^{\beta}$ descends to a very ample line bundle on the coarse moduli space $X$ for some $\alpha,\beta\in\Z_{>0}$. Since the two components of the $\mathbb{G}_m^2$ act on $L_1$ and $L_2$ respectively, we also know that $L_1^{\alpha}$ and $L_2^{\beta}$ are Cartier on $X$. Let $\mathscr{E}=L_1\oplus L_2$ and we may without loss of generality work locally and assume that $B=\Spec R$ is affine where $R$ is a Noetherian ring.\medskip\\
Using the idea of \cite{BOW}, we can embed $\mathcal{X}$ into a toric stack.
\begin{lem} [c.f. 4.2 of \cite{BOW}]
    There is a collection of data $\left(a_0,b_0, \alpha,\beta, V_1, V_2\right)$, where $m_0$ and $\alpha,\beta$ are positive integers and $V_1$ and $V_2$ are finitely generated $R$-modules with $\mathbb{G}_m^2$-action, satisfying\begin{itemize}
        \item [(1)] The algebra $\oplus_{a,b \geq 0} \pi_* (L_1^a\otimes L_2^b)$ is generated by $\oplus_{0 \leq a \leq a_0, 0\leq b\leq b_0} \pi_* (L_1^a\otimes L_2^b) $.
        \item [(2)] The line bundle $L_1^{\alpha}\otimes L_2^{\beta}$ descends to a very ample line bundle on $X$, denoted hereafter by $\OO_X(1)$.
        \item [(3)] $V_1 = H^0\left(X, \OO_X(1)\right)$ is a finitely generated $R$-module and $\mathbb{G}_m^2$-module inducing an immersion
$$
i: X \hookrightarrow \mathbb{P}V_1,
$$
where the $\mathbb{G}_m^2$-action on $\OO_X(1)$ is given by the $\mathbb{G}_m$-actions on $L_1$ and $L_2$.
        \item [(4)] $V_2 = H^0\left(X, \oplus_{0 \leq a\leq a_0,\,0\leq b\leq b_0} \pi_* (L_1^a\otimes L_2^b) \otimes \mathcal{O}_X(1)\right)$ is a finitely generated $R$-module and $\mathbb{G}_m^2$-module inducing a surjection
$$
V_2 \otimes \mathcal{O}_X \rightarrow \oplus_{0 \leq a\leq a_0,\,0\leq b\leq b_0} \pi_* (L_1^a\otimes L_2^b) \otimes \mathcal{O}_X(1).
$$
    \end{itemize}
\end{lem}
\begin{pf}
    First note that (2) is true by definition. For (1), by 3.2 of \cite{BOW}, the algebraic space $$\underline{\Spec}_X\left(\oplus_{a,b \geq 0} \pi_*(L_1^a\otimes L_2^b)\right)$$ is locally of finite type over $X$, so there exist integers $a_0,b_0$ such that $\oplus_{a,b \geq 0} \pi_* (L_1^a\otimes L_2^b)$ is generated by $\oplus_{0 \leq a \leq a_0,0\leq b\leq b_0} \pi_*(L_1^a\otimes L_2^b)$. Since $\OO_X(1)$ is very ample, $\oplus_{0 \leq a \leq a_0,0\leq b \leq b_0} \pi_*(L_1^a\otimes L_2^b) (n)$ is generated by global sections for $n \gg 0$, which we may without loss of generality assume to be true for $n=1$. For (3) and (4), we have $V_1\subset H^0(\OO_X(1))$ and $V_2\subset H^0(\oplus \pi_*(L_1^a\otimes L_2^b)\otimes \OO_X(1))$ by proof of 4.2 in \cite{BOW}, and equality follows from the fact that $\mathcal{X}$ is proper.\qed
\end{pf}
\noindent Now let $\mathbb{H}=\mathbb{A}_{\mathscr{E}}=\underline{\Spec}_{\mathcal{X}}\bigoplus_{n\geq0}S^n\mathscr{E}$ where $\mathscr{E} = L_1\oplus L_2$. Clearly $\mathbb{H}=\mathbb{L}_1\times_{\mathcal{X}}\mathbb{L}_2=\underline{Hom}_{\mathcal{X}}(\mathscr{E},\OO_{\mathcal{X}})$, and hence can be written as \[
    \mathbb{H}=\underline{\Spec}_{\mathcal{X}}\bigoplus_{a,b\geq0}{L}_1^a\otimes{L}_2^b\]
    with coarse moduli space
    \[H=\underline{\Spec}_X\bigoplus_{a,b\geq0}\pi_*({L}_1^a\otimes {L}_2^b).\]
Let
$$\text{mult}^{\alpha,\beta}:\mathbb{H}\to\mathbb{A}_{\OO_X(1)}$$
be the product of $\alpha$-th and $\beta$-th powers of the maps on the two components. Let $\mathbb{I}\subset \mathbb{H}$ be the locus where $\text{mult}^{\alpha,\beta}$ is nonzero, i.e. $\mathbb{I}=\mathcal{P}_{L_1,L_2}$. It is a principal $\mathbb{G}_m^2$-bundle over $\mathcal{X}$. Set $V=V_1\oplus V_2$. We have a map $H\to \mathbb{A}_{V_2}$ defined as follows: let $$\delta:\OO_X(1)\to\bigoplus_{a,b\geq0}\pi_*(L_1^a\otimes L_2^b)$$ be the map defined by $\text{mult}^{\alpha,\beta}$, and then 
$$V_2= H^0\left(\bigoplus_{0\leq a\leq a_0,\,0\leq b\leq b_0}\pi_*(L_1^a\otimes L_2^b)(1)\right)\to \bigoplus_{a,b\geq0}\pi_*(L_1^a\otimes L_2^b)$$ defines $H\to\mathbb{A}_{V_2}$. Let $h:\mathbb{I}\hookrightarrow \mathbb{H}\xrightarrow{\operatorname{CMS}} H$.
\[\begin{tikzcd}
	{\mathbb{I}} & {H\times_X\mathbb{A}_{\OO_X(1)}} & {\mathbb{A}_V}=\mathbb{A}_{V_2}\times \mathbb{A}_{V_1} \\
	& {\underline{\Spec}_X\bigoplus_{a,b,n\geq 0}\pi_*(L_1^{a+n\alpha}\otimes L_2^{b+n\beta})}
	\arrow["{h,\text{mult}^{\alpha,\beta}}", from=1-1, to=1-2]
	\arrow["g", from=1-2, to=1-3]
	\arrow["{=}"{marking, allow upside down}, draw=none, from=1-2, to=2-2]
\end{tikzcd}\]
Let $s\in V_1$ be a section of $\OO(1)$, denote the nonvanishing locus of $s$ in $\mathbb{A}_{V_1}$ by $D(s)$. The corresponding open in $\PP(V_1)$ is denoted $D_+(s)$, and $X_s\coloneqq X\cap D_+(s)$. Look at the Cartesian diagram (denote $\mathbb{A}_{\OO_X(1)}$ by $\mathbb{L}$, $\mathbb{L}^*\coloneqq \mathbb{L}\setminus (0$-section) and $\mathbb{L}^*_s$ means its restriction to $X_s$):
\[\begin{tikzcd}
	{\mathbb{L}_s^*} & {D(s)\subset \mathbb{A}_{V_1}} \\
	{X_s} & {D_+(s)\subset\PP V_1}
	\arrow[from=1-1, to=1-2]
	\arrow[from=1-1, to=2-1]
	\arrow[from=1-2, to=2-2]
	\arrow[hook, from=2-1, to=2-2]
\end{tikzcd}\]
from which we may conclude that $\mathbb{L}_s^*=\mathbb{G}_{m,X_s}$.
\[\begin{tikzcd}
	{\underline{\Spec}_{\mathcal{X}_s}\bigoplus_{a,b\in\Z} L_1^a\otimes L_2^b} & {\mathbb{I}_s} \\
	{\underline{\Spec}_{X_s}\left(\bigoplus_{a,b\geq0}\pi_*(L_1^a\otimes L_2^b)\otimes\bigoplus_{n\in\Z}\pi_*(L_1^{n\alpha}\otimes L_2^{n\beta})\right)} & {H_s\times_{X_s}\mathbb{L}_s^*} \\
	{H^0(\bigoplus_{0\leq a\leq a_0,\,0\leq b\leq b_0}\pi_*(L_1^a\otimes L_2^b)\otimes\OO_X(1))\times \mathbb{L}_s^*} & {\mathbb{A}_{V_2}\times \mathbb{L}_s^*} & {\mathbb{A}_{V_2}\times D(s)} \\
	\\
	& {X_s} & {D_+(s)}
	\arrow["{=}"{description}, draw=none, from=1-1, to=1-2]
	\arrow["{h,\text{mult}^{\alpha,\beta}}", hook, from=1-2, to=2-2]
	\arrow["{=}"{marking, allow upside down}, draw=none, from=2-2, to=2-1]
	\arrow["j", hook, from=2-2, to=3-2]
	\arrow["{=}"{marking, allow upside down}, draw=none, from=3-2, to=3-1]
	\arrow[from=3-2, to=3-3]
	\arrow[from=3-2, to=5-2]
	\arrow["\lrcorner"{anchor=center, pos=0.125}, draw=none, from=3-2, to=5-3]
	\arrow[from=3-3, to=5-3]
	\arrow["i", hook, from=5-2, to=5-3]
\end{tikzcd}\]
\begin{cl} [c.f. \cite{BOW} 4.5]
    The maps $(h,\text{mult}^{\alpha,\beta})$ and $j$ are both closed immersions.
\end{cl}
\begin{pf}
We first show that $j:H_s\times_{X_s}\mathbb{L}_s^*\to \mathbb{A}_{V_2}\times \mathbb{L}_s^*$ is a closed immersion. Clearly $\mathbb{A}_{V_2}\times \mathbb{L}_s^*$ and $H_s\times_{X_s}\mathbb{L}_s^*$ are both affine over $X_s$ and $j$ is given by the sheaf map
$$\gamma:\OO_{X_s}[t^{\pm}]\otimes_R S^{\bullet}V_2\to\OO_{X_s}[t^{\pm}]\otimes_{\OO_{X_s}} \bigoplus_{a,b\geq0}\pi_*(L_1^a\otimes L_2^b)|_{X_s}$$
where the right hand side can be identified as
$$\bigoplus_{n\in\Z \text{ and } a,b\geq 0}\left(\OO_{X_s}(n)\otimes_{\OO_{X_s}}\pi_*(L_1^a\otimes L_2^b)|_{X_s}\right)$$ sending $t$ to the nowhere vanishing section $s$ in $\OO_{X_s}(1)$. Under this identification, $\gamma$ is defined by sending $v\in V_2$ to its image in the $n=1$ summand. By (1) in the previous lemma, we know that the right hand side is generated by $\bigoplus_{0\leq a\leq a_0,\,0\leq b\leq b_0}\pi_*(L_1^a\otimes L_2^b)$ as $\OO_{X_s}[t^{\pm}]$-algebra, and hence due to (4) it's generated by the image of left hand side. So $j$ is a closed immersion.\medskip\\
Next, we show that $(h,\text{mult}^{\alpha,\beta})$ is a closed immersion. Note that since $\mathbb{I}=\mathcal{P}_{L_1,L_2}$ is already an algebraic space, the coarse moduli map doesn't do anything on it, so $h$ is an open immersion. Therefore we have the following Cartesian diagram:
\[\begin{tikzcd}
	{\mathbb{I}_s} && {H_s\times_{X_s}\mathbb{L}_s^*} \\
	{H_s} && {H_s\times_{X_s} \mathbb{L}_s}
	\arrow["{h,\operatorname{mult}^{\alpha,\beta}}", from=1-1, to=1-3]
	\arrow["h"', from=1-1, to=2-1]
	\arrow[from=1-3, to=2-3]
	\arrow["{\operatorname{id},\operatorname{mult}^{\alpha,\beta}}"', from=2-1, to=2-3]
\end{tikzcd}\]
since the bottom horizontal map is a closed immersion, so is $(h,\text{mult}^{\alpha,\beta})$.\qed
\end{pf}
\noindent Now we have an immersion $\mathbb{I}_s\hookrightarrow \mathbb{A}_{V_2}\times D(s)$, and so $\mathbb{I}\to \mathbb{A}_V$ is an immersion. Let $x\in I$ be a geometric point, and $s\in V_1$ be a section not passing through $x$. Let $(D(s)\times \mathbb{A}_{V_2})_{k(x)}$ be the fiber over $\Spec \,k(x)\to\Spec R$. We would like to show that the orbit of $x$ in $(D(s)\times \mathbb{A}_{V_2})_{k(x)}$ is closed. Let $y\in X_s$ be the image of $x$. Since $\PP V_1$ is separated over $R$, this is equivalent to show orbit of $x$ is closed in $(D(s)\times \mathbb{A}_{V_2})_{y}$, i.e. closed in $(\mathbb{A}_{V_2}\times \mathbb{L}_s^*)_y$, i.e. closed in $\mathbb{I}_y$, but $\mathbb{I}$ over $\mathcal{X}$ is a $\mathbb{G}_m^2$-torsor so orbit is $\mathbb{I}_y$. Also since diagonal of $\mathcal{X}$ is finite, stabilizer group of $x$ is finite, so $x\in\mathbb{A}_V$ is stable. Thus, for the twisted affine GIT problem:\begin{itemize}
    \item [(1)] action is given by $\mathbb{G}_m^2\acts \mathbb{A}_V:((z_1,z_2),v) \mapsto \sum z_1^az_2^b v_{a,b}$  where $v = \sum v_{a,b}$ is a decomposition of $v$ such that $v_{a,b}\in H^0(\pi_*(L_1^a\otimes L_2^b))$
    \item [(2)] character is given by $\text{mult}^{\alpha,\beta}:(z_1,z_2)\mapsto z_1^{\alpha}z_2^{\beta}$
    \item [(3)] $\QP(V)$ denotes the stack quotient of the stable locus by $\mathbb{G}_m^2$.
\end{itemize}
\begin{thm}
    We have an immersion $\mathbb{I}\hookrightarrow \mathbb{A}_V^{s,\,\operatorname{mult}^{\alpha,\beta}}$. Taking $\mathbb{G}_m^2$ quotient we get $\mathcal{X}\hookrightarrow \QP(V)$.
\end{thm}
\begin{rem}
    Recall that $\QP(V)$ is proper only if the stable locus $\mathbb{A}_V^{s}$ in $\mathbb{A}_V$ equals the semistable locus $\mathbb{A}_V^{ss}$ since $[\mathbb{A}_V^{ss,\chi}/\mathbb{G}_m^2]$ is a proper tame stack where the character $\chi$ is given by $\operatorname{mult}^{\alpha,\beta}$ and $[\mathbb{A}_V^{s,\chi}/\mathbb{G}_m^2]\subset [\mathbb{A}_V^{ss,\chi}/\mathbb{G}_m^2]$ is open. But by the Hilbert--Mumford criterion for twisted affine GIT (Proposition 2.5 of \cite{king}), we may take some $x\in V$, write it as $x=x_1+x_2$ where $x_1\in V_1$ and $x_2\in V_2$, such that $x_1\neq 0$ and $x_2 $ has a nonzero component in $H^0(\pi_*(L_1^a\otimes L_2^b))$ only when $a\beta-b\alpha\geq 0$. Such $x$ will be a strictly semistable point in $\mathbb{A}_V$ by looking at the 1 parameter subgroup $\G_m \to \G_m^2: t\mapsto (t^{\beta},t^{-\alpha})$.
\end{rem}
\begin{rem}
    [Kirwan blowups] $\mathcal{X}$ does embed into a tame stack with proper good moduli space. Let $\rho:\widetilde{\mathbb{A}}_V\to \mathbb{A}_V$ be the Kirwan blowup with action by $\mathbb{G}_m^2$. On $\widetilde{\mathbb{A}}_V$, the semistable locus of the $\mathbb{G}_m^2$-action coincides with the stable locus, i.e. $\widetilde{\mathbb{A}}_V^s = \widetilde{\mathbb{A}}_V^{ss}$. There is an open embedding of $\mathbb{A}_V^s\hookrightarrow \widetilde{\mathbb{A}}_V^s$, and thus there's a closed embedding $\mathcal{X}\hookrightarrow [\widetilde{\mathbb{A}}_V^{s}/\mathbb{G}_m^2]$ into a proper stack. See \cite{kirwanbu}.
    \[\begin{tikzcd}
	{\mathcal{X}=[\mathbb{I}/\mathbb{G}_m^2]} & {[\mathbb{A}_V^s/\mathbb{G}_m^2]} & {[\mathbb{A}_V^{ss}/\mathbb{G}_m^2]} \\
	& {[\widetilde{\mathbb{A}}_V^{s}/\mathbb{G}_m^2]} & {[\widetilde{\mathbb{A}}_V^{ss}/\mathbb{G}_m^2]}
	\arrow[hook, from=1-1, to=1-2]
	\arrow["\shortmid"{marking}, hook, from=1-1, to=2-2]
	\arrow[hook, from=1-2, to=1-3]
	\arrow["\circ"{marking, allow upside down}, hook', from=1-2, to=2-2]
	\arrow["{=}"', from=2-2, to=2-3]
	\arrow[from=2-3, to=1-3]
\end{tikzcd}\]
\end{rem}
This argument generalizes to $k$-fold. Assume that $f:\mathcal{X}\to B$ is a $k$-cyclotomic stack uniformized by $L_1,...,L_k$. $(\alpha_1,...,\alpha_k)\in\Z^d_{>0}$ is such that $L_1^{\alpha_1}\otimes\cdots\otimes L_k^{\alpha_k}$ descends to a relatively very ample line bundle $\OO_X(1)$ on the coarse space family $X\to B$. Then, under the same twisted affine GIT problem with\begin{itemize}
    \item $V = V_1 \oplus V_2$, where $V_1=H^0(X,\OO_X(1))$ and $V_2 = H^0(X, \bigoplus_{0\leq a_i\leq a_i'} \pi_*(L_1^{a_1}\otimes\cdots\otimes L_k^{a_k})\otimes \OO_X(1))$ for some positive integers $a_1',...,a_k'$;
    \item character $\text{mult}^{\underline{\alpha}}: \G_m^k\to\G_m:(z_1,...,z_k)\mapsto z_1^{\alpha_1}z_2^{\alpha_2}\cdots z_k^{\alpha_k}$
\end{itemize} we analogously conclude:
\begin{thm}
    Let $f:\mathcal{X}\to B$ be a proper $k$-cyclotomic stack, uniformized by $L_1,...,L_k$ (the morphism $\mathcal{X}\to B\mathbb{G}_m^k$ is representable) with coarse space map $\pi:\mathcal{X}\to X$ and coarse space family $\bar{f}:X\to B$. If $L_1^{\otimes\alpha_1}\otimes L_2^{\otimes\alpha_2}\otimes\cdots\otimes L_k^{\otimes\alpha_k}$ descends to a very ample line bundle on $X$ (the $i$-th component of $\mathbb{G}_m^k$ acts on $L_i$ only), then we have closed immersion $\mathcal{X}\hookrightarrow \QP(V)$.
    \label{embedding}
\end{thm}
\section{Construction of Moduli Space}
\noindent For any family of $k$-uniformized $k$-cyclotomic stack $(f:\mathcal{X}\to B,L_1,...,L_k)$, if there is a a direction $(\alpha_1,...,\alpha_k)\in\Z^k$ such that $L_1^{\alpha_1}\otimes\cdots\otimes L_k^{\alpha_k}$ descends to a relatively ample line bundle on $X$ over $B$, then we could construct the corresponding moduli space. This section gives the construction of the stack $\text{Sta}^{L_1,\,...,\, L_k,\,\underline{\alpha}}$, which is in most cases the first step in the moduli construction. One could modify the stack by taking unions, rigidifying, picking out sub-loci satisfying certain conditions to meet their specific needs. For example, in \cite{bp}, we'll use this section to construct the moduli space of Alexeev--Filipazzi--Inchiostro (AFI) stable families \cite{AFI} and then prove that $\mathcal{M}^{\text{AFI}}_{n,p,I}$ is a DM stack of finite type, and has projective coarse moduli space. For clarity, we continue to assume that $k=2$ in this section but the general case has no difference. This is a generalization of Section 3 in \cite{AH}.\medskip\\ 
In this section, we assume that $L_1^{\alpha_1}\otimes\cdots\otimes L_k^{\alpha_k}$ (or when $k=2$, $L_1^{\alpha_1}\otimes L_2^{\alpha_2}$) descends to a relatively very ample line bundle on the coarse space family $X\to B$. All GIT problems are twisted by the character $\text{mult}^{\underline{\alpha}}:\G_m^k\to \G_m:(z_1,...,z_k)\mapsto z_1^{\alpha_1}... z_k^{\alpha_k}$ as in \ref{embedding}. We'll omit the character when writing down stable or semistable locus, e.g. $\mathbb{A}_V^{s}\coloneqq\mathbb{A}_V^{s,\,\text{mult}^{\underline{\alpha}}}$.
\begin{defn}
    Let $B$ be a scheme of finite type and $V$ a trivial locally-free $\mathcal{O}_{B}$-module. Suppose that $\mathbb{G}_m^2$ acts $\mathcal{O}_B$-linearly on $V$ with positive weights $\left((v_1,w_1), \cdots, (v_m,w_m)\right)$. We have a decomposition of $\mathcal{O}_B$-modules
$$V \simeq V_{v_1,w_1} \oplus \ldots \oplus V_{v_m,w_m}$$
with
$$\rho(a,b)(x)=a^{v_i}b^{w_i} x, \quad \forall\,x \in V_{v_i,w_i},\, (a,b)\in\mathbb{G}_m^2$$
Let $\mathbb{A}$ be the associated vector bundle $\underline{\Spec}_B \sym_B^{\bullet}V$, $\rho:\mathbb{G}_m^2\times_B\mathbb{A}\to\mathbb{A}$. We define \[\mathcal{P}(\rho)\coloneqq [\mathbb{A}^{ss}/\mathbb{G}_m^2]\] Since we are taking quotient on the semi-stable locus, $\mathcal{P}(\rho)$ has a proper good moduli space. Define $\QP(\rho)$ to be $[\mathbb{A}^s/\mathbb{G}_m^2]$. $\QP(\rho)\subset \mathcal{P}(\rho)$ is an open substack. Denote by $\OO(1,0)$ and $\OO(0,1)$ the two tautological line bundles given by the map $[\mathbb{A}^{ss}/\mathbb{G}_m^2]\to B\mathbb{G}_m^2$.
\end{defn}
\begin{defn}
    Define \[\hilb_{\mathcal{Q}\mathcal{P}(\rho)}(T)\coloneqq \{\text{closed substacks }\mathcal{X}\subset \QP(\rho)_T \text{ with }\mathcal{X}\text{ flat and proper over }T\}\] For a function $\mathfrak{F}:\Z^2\to\Z$, if we further impose Hilbert function on the fibers $\chi(\mathcal{X}_t,\OO_{\mathcal{X}_t}(a,b))=\mathfrak{F}(a,b)$ for every $t\in T$, then we get $\hilb_{\QP(\rho),\,\mathfrak{F}}$.
\end{defn}
\begin{defn}
    We define Sta$^{L_1,L_2,\,c}$ to be a category fibered in groupoids over the category of schemes, and\begin{itemize}
        \item [(i)] an object of Sta$^{L_1,L_2,\,c}(B)$ is a proper family of bicyclotomic stacks $\mathcal{X}$, uniformized by $L_1 $ and $L_2$, with $ L_1+c L_2$ being a polarizing direction (a positive multiple of $ L_1+c L_2$ descends to a relatively ample line bundle on the coarse moduli $X$ over $B$).
        \item[(ii)] an arrow from $(\mathcal{X},L_1,L_2)$ to $(\mathcal{X}',L_1',L_2')$ is a fibered diagram
        \[\begin{tikzcd}
	{\mathcal{X}} & {\mathcal{X}'} \\
	B & {B'}
	\arrow["\phi", from=1-1, to=1-2]
	\arrow[from=1-1, to=2-1]
	\arrow["\lrcorner"{anchor=center, pos=0.125}, draw=none, from=1-1, to=2-2]
	\arrow[from=1-2, to=2-2]
	\arrow[from=2-1, to=2-2]
\end{tikzcd}\]
and isomorphisms $\gamma_1:L_1\to \phi^*L_1'$, $\gamma_2:L_2\to \phi^*L_2'$.
    \end{itemize}
\end{defn}
\begin{defn} \label{normalembedding}
    A closed substack $\mathcal{X}\subset\QP(\rho)$ is said to be $(\alpha,\beta)$-normally embedded (where $\beta/\alpha=c$) if \begin{itemize}
        \item[(1)]\indent for any $i,j$, the homomorphism $V_{v_i,w_i}\to H^0(\mathcal{X},\OO_{\mathcal{X}}(v_i,w_i))$ is an isomorphism
        \item[(2)] if $(a,b)\geq (\alpha,\beta)$, then $H^j(\mathcal{X},\OO_{\mathcal{X}}(a,b))=0$ for all $j\geq 1$
        \item[(3)] there is a map $\sym V\to\bigoplus_{a,b\geq 0}H^0(\mathcal{X},\OO_{\mathcal{X}}(a,b))$. It's a surjection onto the components where $(a,b)\geq (\alpha,\beta)$
        \item[(4)] the line bundle $\OO_{\mathcal{X}}(\alpha,\beta)$ descends to a very ample line bundle on $X$.
    \end{itemize}
\end{defn}
\begin{rem}
    We wanted to say that (2) and (3) hold for $(a,b)\geq (v_i,w_i)$ for some pair of weights indexed by $i$, but in our representation that's equivalent to $(a,b)\geq(\alpha,\beta)$. Notice also that Condition (2) implies that $H^j(X,\pi_*\OO(a,b))=0$ for $j>0$ from the Leray spectral sequence.
\end{rem}
\begin{defn}
    A stack $\mathcal{X}$ uniformized by $L_1,L_2$ is $(\alpha,\beta)$-normally embeddable in $\QP(\rho)$ if
there exists a $(\alpha,\beta)$-normal embedding $\mathcal{X}\subset\QP(\rho)$ such that $\OO_{\mathcal{X}} (a,b) =L_1^a\otimes L_2^b$ for all $a,b\in\Z$.
\end{defn}

\begin{prop}
    Let $(f:\mathcal{X}\to B,L_1,L_2)$ be a family of 2-uniformized bicyclotomic stacks over a Noetherian base $B$ and $L_1^{\alpha_0}\otimes L_2^{\beta_0}$ descends to a very ample line bundle on $X$. Then $f:\mathcal{X}\to B$ is normally embeddable.
\end{prop}
\begin{pf}
     We assumed that $L_1^{\alpha_0}\otimes L_2^{\beta_0}$ is a polarizing direction. Since the relative ample cone for $X\to B$ in the plane $\Q\langle L_1^{N_1},L_2^{N_2}\rangle$ is open, we could find a set of basis $\{L_1',L_2'\}$ for the lattice $\Z\langle L_1,L_2\rangle$ where $L_1',L_2'$ sit in the relative ample cone. $f:\mathcal{X}\to B$ is still uniformized by $L_1' $ and $L_2'$, so we could simply assume $L_1$ and $L_2$ are polarizing to start with. Now choose $(\alpha,\beta)\in\Z^2$ be such that (1) $L_1^{\alpha}\otimes L_2^{\beta}$ descends to the coarse space as a relatively very ample line bundle and (2) for any $(a,b)\geq (\alpha,\beta)$ we have $H^i(\pi_*(L_1^a\otimes L_2^b))=0$. The second condition is possible due to 1.2.22 of \cite{larz1}, by running the argument with $\mathcal{F}\coloneqq \pi_*L_1^{r}$ and $\mathcal{G}\coloneqq \pi_*L_2^s$ for each $0\leq r<N_1$ and $0\leq s<N_2$, and then take the maximum of all.\qed
\end{pf}

\begin{prop}
    Fix a function $\mathfrak{F}:\Z^2\to\Z$ and positive integers $a_0$ and $b_0$. For $\alpha\leq a\leq \alpha+a_0$ and $\beta\leq b\leq \beta+b_0$, let $V_{a,b}$ be a free module with rank $\mathfrak{F}(a,b)$. Consider the $\mathbb{G}_m^2$-representation $\rho=\rho_{\alpha,\beta,a_0,b_0}$ on
    $$V=V_{\alpha,\beta}\oplus\bigoplus_{\alpha\leq a\leq \alpha+a_0,\,\beta\leq b\leq \beta+b_0}V_{a,b}$$
    \begin{itemize}
        \item [(1)] there exists an open subscheme $H\subset\hilb_{\QP(\rho),\,\mathfrak{F}}$ parametrizing $(\alpha,\beta)$-normally embedded substacks;
        
        \item [(2)] the subcategory Sta$^{L_1,L_2,\,c}_{\mathfrak{F}}(\rho)$ of objects with fibers $(\alpha,\beta)$-normally embeddable in $\QP(\rho)$ and $\chi(\mathcal{X},L_1^a\otimes L_2^b)=\mathfrak{F}(a,b)$ satisfies Sta$^{L_1,L_2,\,c}_{\mathfrak{F}}(\rho)=[H/G]$ where $G=\text{Aut}_{\rho}(V)$.
    \end{itemize}
\end{prop}
\begin{rem}
    Note that due to condition (2) in the definition of $(\alpha,\beta)$-normal embedding, $\chi(\mathcal{X},L_1^a\otimes L_2^b)$ is equal to $h^0(\mathcal{X},L_1^a\otimes L_2^b)$ for $\alpha\leq a\leq \alpha+a_0,\beta\leq b\leq \beta+b_0$.
\end{rem}
\begin{pf} [of the proposition]
    For (1), all the conditions in the definition of $(\alpha,\beta)$-normal embedding are open conditions due to the upper semicontinuity theorem and the theorem on cohomology and base change. Indeed, if one fiber $X_t$ over $t\in B$ has vanishing higher cohomologies for $\OO_X(a,b)|_{X_t}$ whenever $(a,b)\geq (\alpha,\beta)$, this must also be true in a neighbourhood of $t\in B$. Then the pushforward $f_*(\OO_{\mathcal{X}}(a,b))$ will be locally free near $t\in B$, and in an open neighbourhood $V\subset B$ containing $t$, $X_V\to \PP(f_*(\OO_{\mathcal{X}}(\alpha,\beta)|_V)^{\vee}$ is a closed embedding. \medskip\\Next we prove (2). The map $B\to[H/G]$ is the same data as of a map 
    $$\mathcal{Y}\to Q\to B$$ where $Q$ is a principal $G$-bundle over $B$ and $\mathcal{Y}\hookrightarrow \QP(\rho)\times Q$ a family of $(\alpha,\beta)$-normally embedded substacks. There is a principal $\mathbb{G}_m^2$-bundle 
    $$\mathcal{P}_{\OO_{\mathcal{Y}}(1,0),\OO_{\mathcal{Y}}(0,1)}\to\mathcal{Y}$$
    The $G$-action lifts to this $\mathbb{G}_m^2$-bundle and commutes with the $\mathbb{G}_m^2$-action. The free quotient $\mathcal{P}_{\OO_{\mathcal{Y}}(1,0),\OO_{\mathcal{Y}} (0,1)}/G$ retains a $\mathbb{G}_m^2$-action with finite stabilizers, and we denote the quotient $\mathcal{X}$. Let\begin{itemize}
        \item [(1)] $L_1$ be the line bundle associated to $\mathcal{P}_{\OO_{\mathcal{Y}}(1,0)}/G\to \mathcal{X}$
        \item [(2)] $L_2$ be the line bundle associated to $\mathcal{P}_{\OO_{\mathcal{Y}}(0,1)}/G\to \mathcal{X}$
    \end{itemize}
    Here we identify $\mathcal{X}=\mathcal{Y}/G$, $L_1=\mathbb{A}_{\OO_{\mathcal{Y}}(1,0)}/G$ and $L_2=\mathbb{A}_{\OO_{\mathcal{Y}}(0,1)}/G$. In summary, we have:
    \[\hspace{-1cm}\begin{tikzcd}
	{\mathcal{P}_{\OO_{\mathcal{Y}}(1,0),\OO_{\mathcal{Y}}(0,1)}} & {\mathbb{A}_{\OO_{\mathcal{Y}}(0,1)}\times \mathbb{A}_{\OO_{\mathcal{Y}}(1,0)}} & {\mathbb{A}_{\OO(1,0)}\times\mathbb{A}_{\OO(0,1)}\times Q} & {\mathbb{A}_{\OO(1,0)}\times\mathbb{A}_{\OO(0,1)}\times H} \\
	{\mathcal{P}_{L_1.L_2}} & {\mathcal{Y}} & {\QP(\rho)\times Q} & {\QP(\rho)\times H} \\
	& {\mathcal{X}} & Q & H \\
	&& B & {[H/G]}
	\arrow[hook, from=1-1, to=1-2]
	\arrow[from=1-1, to=2-1]
	\arrow["{\mathbb{G}_m^2}"', from=1-1, to=2-2]
	\arrow[hook, from=1-2, to=1-3]
	\arrow[from=1-2, to=2-2]
	\arrow[from=1-3, to=1-4]
	\arrow[from=1-3, to=2-3]
	\arrow[from=1-4, to=2-4]
	\arrow["{\mathbb{G}_m^2}"', from=2-1, to=3-2]
	\arrow[hook, from=2-2, to=2-3]
	\arrow["G", from=2-2, to=3-2]
	\arrow[from=2-2, to=3-3]
	\arrow[from=2-3, to=2-4]
	\arrow[from=2-3, to=3-3]
	\arrow[from=2-4, to=3-4]
	\arrow[from=3-2, to=4-3]
	\arrow[from=3-3, to=3-4]
	\arrow["G"', from=3-3, to=4-3]
	\arrow[from=3-4, to=4-4]
	\arrow[from=4-3, to=4-4]
\end{tikzcd}\]
and in this diagram all parallelograms are Cartesian. Hence the geometric fibers of $\mathcal{X}\to B$ are the same as those for $\mathcal{Y}\to Q$, they are bicyclotomic, uniformized by $\OO(1,0)$ and $\OO(0,1)$ (due to commutativity of $G$- and $\mathbb{G}_m^2$-actions), and $(\alpha,\beta)$-normally embedded in $\QP(\rho)$. The $(\alpha,\beta)$-normal embedding condition also implies that $\OO_{\mathcal{X}}(\alpha,\beta)$ descends to a relatively very ample line bundle on the coarse moduli $X\to B$.\medskip\\
Since the construction works for arrows $B'\to B\to [H/G]$, we get a functor 
$$[H/G]\to\operatorname{Sta}^{L_1,L_2,\,c}_{\mathfrak{F}}(\rho).$$
Conversely we start from a family $\mathcal{X}\to B$ in $\operatorname{Sta}^{L_1,L_2,\,c}_{\mathfrak{F}}(\rho).$ Consider the locally free sheaf (since all higher cohomologies vanish) 
$$W=f_*(L_1^{\alpha}\otimes L_2^{\beta})\oplus\bigoplus_{\alpha\leq a\leq \alpha+a_0,\beta\leq b\leq \beta+b_0}f_*(L_1^a\otimes L_2^b)$$ 
Each fiber of $W$ is isomorphic to $V$ with representation $\rho$ as $\mathbb{G}_m^2$-space. By results in the previous section we get embedding $\mathcal{X}\hookrightarrow \QP(W)$. The principal $G$-bundle $Q=\text{Isom}_{\mathbb{G}_m^2}(W,V_B)$ satisfies $W\times_B Q\cong V_Q$ and hence $\QP(W)\times_BQ\cong\QP(\rho)\times Q$. Write $\mathcal{Y}=\mathcal{X}\times_BQ\hookrightarrow\QP(\rho) \times Q$, we get an $(\alpha,\beta)$-normal embedding by definition, and this is $G$-equivariant, so we get an object in $[H/G]$. Again this construction works for a map $B'\to B\to \operatorname{Sta}^{L_1,L_2,\,c}_{\mathfrak{F}}(\rho)$, hence we obtain a functor in the reverse direction $\operatorname{Sta}^{L_1,L_2,\,c}_{\mathfrak{F}}(\rho)\to[H/G]$.\qed
\end{pf}

\begin{thm}
        $\operatorname{Sta}^{L_1,L_2,\,c}_{\mathfrak{F}}$ is algebraic and locally of finite type.
\end{thm}
\begin{pf}
    We first show that $\operatorname{Sta}^{L_1,L_2,\,c}_{\mathfrak{F}}(\rho_{\alpha,\beta,a_0,b_0})\subset\operatorname{Sta}^{L_1,L_2,\,c}_{\mathfrak{F}}$ is an open inclusion. We have argued before that the vanishing of higher cohomologies and very ample-ness are both open conditions. Note also that in any family $\mathcal{X}\to B$, if $a_0,b_0\in\Z_{>0}^2$ is such that some fiber $\mathcal{X}_b\hookrightarrow \QP(\rho_{\alpha,\beta,a_0,b_0})$, then since the locus in $B$ where $f_*(L_1^{\alpha}\otimes L_2^{\beta})\oplus \bigoplus_{\alpha\leq a\leq\alpha+a_0, \beta\leq b\leq \beta+b_0}f_*(L_1^a\otimes L_2^b)$ is isomorphic to $V_B$ where $V$ is the space $V_{\alpha,\beta}\oplus \bigoplus_{\alpha\leq a\leq\alpha+a_0, \beta\leq b\leq \beta+b_0}V_{a,b}$ is open ($\mathfrak{F}$ controls the dimension of each piece due to vanishing of higher cohomologies). Therefore, the same choices for $\alpha,\beta,a_0,b_0$ works for an entire open neighbourhood of $b$ in the base $B$. Since $\operatorname{Sta}^{L_1,L_2,\,c}_{\mathfrak{F}}=\bigcup \operatorname{Sta}^{L_1,L_2,\,c}_{\mathfrak{F}}(\rho_{\alpha,\beta,a_0,b_0})$ is a union of algebraic stacks locally of finite type, $\operatorname{Sta}^{L_1,L_2,\,c}_{\mathfrak{F}}$ is algebraic and locally of finite type. \qed
\end{pf}
In general, we have:
\begin{thm}
    $\operatorname{Sta}^{L_1,\,...,\,L_k,\,\underline{\alpha}}_{\mathfrak{F}}$ is algebraic and locally of finite type. Here $\underline{\alpha}$ means that $L_1^{\alpha_1}\otimes \cdots\otimes L_k^{\alpha_k}$ descends to a relatively ample line bundle on the coarse space. 
    \label{stastack}
\end{thm}

\section{Koll\'ar Family of Seifert $\mathbb{G}_m^d$-Bundles}
\noindent In this section we generalize the ordinary mAH construction to include the case of $L^{[n]}(\lfloor nD\rfloor)$ where $L$ is a reflexive rank 1 sheaf and $D$ is a boundary divisor on $X$.
\begin{defn}
    By a Koll\'ar family of Seifert $\mathbb{G}_m^d$-bundles, we mean\begin{itemize}
        \item $f:X\to B$ a flat family of equi-dimensional (connected) reduced S2 schemes;
        \item $\{L_{\underline{n}}\}_{\underline{n}\in\Z^d}$ a lattice of coherent sheaves on $X$ flat over $B$ such that\begin{itemize}
            \item [(1)] for any $b\in B$ and $\underline{n}\in\Z^d$, $L_{\underline{n}}|_{X_b}$ is a reflexive sheaf of rank 1, and $L_{\underline{0}}=\OO_X$.
            \item [(2)] there exists an open subset $U\subset X$ such that for all $\underline{n}\in\Z^d$, $L_{\underline{n}}$ is invertible on $U$.
            \item [(3)] there exist compatible multiplication maps $$L_{\underline{n}}\otimes L_{\underline{m}}\to L_{\underline{m}+\underline{n}}$$ for all $\underline{m},\underline{n}\in\Z^d$. There exists an open subset $V\subset U$ where all multiplication maps are isomorphisms. For arbitrary base change $g:B'\to B$,
            \[\begin{tikzcd}
	{X'} & X \\
	{B'} & B
	\arrow["\varphi", from=1-1, to=1-2]
	\arrow["{f'}", from=1-1, to=2-1]
	\arrow["\lrcorner"{anchor=center, pos=0.125}, draw=none, from=1-1, to=2-2]
	\arrow["f", from=1-2, to=2-2]
	\arrow["g", from=2-1, to=2-2]
\end{tikzcd}\]
all the pullbacks $\varphi^*L_{\underline{n}}$ are reflexive on the fibers of $X'$, and 
$$\varphi^*\Ln \otimes\varphi^*\Lm\to\varphi^*\Lmn$$ is the multiplication structure on $X'\to B'$.
\item [(4)] for any $b\in B$, there exists $M_1,...,M_d\in\Z_{>0}$ such that $L_{M_i\underline{e}_i}|_{X_b}$ is invertible for every $i$. The multiplication maps \[\Ln |_{X_b}\otimes L_{M_i\underline{e}_i}|_{X_b}\to L_{M_i\underline{e}_i+\underline{n}}|_{X_b}\] are isomorphisms for all $i,\underline{n}$.
        \end{itemize}
    \end{itemize}
    
\end{defn}

\begin{rem}
    As before, if $L_{M_i\underline{e}_i}|_{X_b}$ is locally free, then it's locally free in a neighbourhood of $b\in B$. So in particular, if we assume that $B$ is Noetherian, then we could choose the $M_i$'s that work for the entire family. For those $M_i$'s, \[\text{mult}:\Ln \otimes L_{M_i\underline{e}_i}\cong L_{\underline{n}+M_i\underline{e}_i}\] are isomorphisms for any $\underline{n}\in\Z^d$.
\end{rem}

\noindent A morphism $(\varphi,\alpha_{\underline{n}})$ from a Koll\'ar family of Seifert $\mathbb{G}_m^d$-bundles $(X'\to B',\{\Ln'\}_{\underline{n}\in\Z^d})$ to $(X\to B,\{\Ln\}_{\underline{n}\in\Z^d})$ is given by a Cartesian square
\[\begin{tikzcd}
	{X'} & X \\
	{B'} & B
	\arrow["\varphi", from=1-1, to=1-2]
	\arrow["{f'}", from=1-1, to=2-1]
	\arrow["\lrcorner"{anchor=center, pos=0.125}, draw=none, from=1-1, to=2-2]
	\arrow["f", from=1-2, to=2-2]
	\arrow["g", from=2-1, to=2-2]
\end{tikzcd}\]
and isomorphisms $\alpha_{\underline{n}}: \Ln'\to \varphi^*\Ln$ for any $\underline{n}\in\Z^d$ that are compatible with the multiplication maps, i.e. the diagram
\[\begin{tikzcd}
	& {\Ln' \otimes \Lm'} \\
	{\Lmn'} & {\varphi^*\Ln \otimes\varphi^*\Lm} \\
	& {\varphi^*\Lmn}
	\arrow["{\text{mult}}"', from=1-2, to=2-1]
	\arrow["{\alpha_{\underline{n}}\otimes\alpha_{\underline{m}}}", from=1-2, to=2-2]
	\arrow["{\alpha_{\underline{n}+\underline{m}}}"', from=2-1, to=3-2]
	\arrow["{\varphi^*\text{mult}}", from=2-2, to=3-2]
\end{tikzcd}\] commutes. The objects and arrows define the Koll\'ar category of $d$ $\Q$-line bundle sequences.
\begin{rem}
    Whether to impose connectedness condition in fiber depends on whether connectedness condition was included in the definition of orbispace because later we'll establish the equivalence between the two categories.
\end{rem}
On the Koll\'ar family, $\bigoplus_{n\in\Z^d}\Ln$ forms a $\Z^d$-graded algebra under the multiplication structure. As before, we perform the Abramovich--Hassett construction:
\begin{defn}
    The $\mathbb{G}_m^d$-space of $\Ln$ is defined as
    $$\mathcal{P}_{\{\Ln\}}\coloneqq \underline{\Spec}_X\left(\oplus_{\underline{n}\in\Z^d} \Ln\right)$$
    where the $\mathbb{G}_m^d$-action on $\mathcal{P}_{\{\Ln\}}$ is given by the $\Z^d$-grading. The Abramovich--Hassett stack of $\{\Ln\}$ on $X$ is defined to be its stack $\mathbb{G}_m^d$-quotient
    $$\mathcal{X}_{\{\Ln\}}\coloneqq AH(X,\{\Ln\})\coloneqq \left[\mathcal{P}_{\{\Ln\}}\middle/\mathbb{G}_m^d\right]$$ This principal $\mathbb{G}_m^d$-bundle $\mathcal{P}_{\{\Ln\}}\to \mathcal{X}_{\{\Ln\}} $ is naturally associated with line bundles $\mathcal{L}_1,\cdots, \mathcal{L}_d$, with $\mathcal{L}_i$ coming from the $i$-th component of $\mathbb{G}_m^d$, i.e. it corresponds to the direction $\underline{e}_i\in\Z^d$.
\end{defn}
\begin{prop} [Functor] We have the following properties of $\mathcal{X}_{\{\Ln\}}$: \begin{itemize} 
\item [(1)] The family $\mathcal{X}_{\{\Ln\}}\to B$ is a family of $d$-cyclotomic orbispaces with S2 fibers.
\item[(2)] $\pi:\mathcal{X}_{\{\Ln\}}\to X$ is the coarse moduli space map.
\item [(3)] For any $\underline{n}=(n_1,\cdots,n_d)\in\Z^d$, 
$$\Ln = \pi_*(\mathcal{L}_1^{n_1}\otimes\cdots\otimes \mathcal{L}_d^{n_d})$$
\item [(4)] Given a morphism in the Koll\'ar category $(\varphi,\alpha_{\underline{n}})$ from $(X'\to B',\{\Ln'\})$ to $(X\to B,\{\Ln\})$, we get canonically isomorphisms $\mathcal{P}_{\{\Ln'\}}\cong\varphi^*\mathcal{P}_{\{\Ln\}}$ and $\mathcal{X}_{\{\Ln'\}}\cong\varphi^* \mathcal{X}_{\{\Ln\}}$.
\end{itemize}
\end{prop}
\begin{pf}
(1) Locally over $B$, we know that $L_{M_i\underline{e}_i}$ is locally free on $X$ for every $1\leq i\leq d$. For any $\underline{n}=(n_1,\cdots,n_d)\in\Z^d$, the multiplication maps 
$$\text{mult}:(L_{M_1\underline{e}_1})^{\otimes n_1}\otimes\cdots\otimes (L_{M_d\underline{e}_d})^{\otimes n_d}\to L_{(M_1n_1,\, ...,\,M_dn_d)}$$
are all isomorphisms by assumption. Due to $\Z^d$-grading, we have $\mathbb{G}_m^d$-equivariant maps
\[\mathcal{P}_{\{\Ln\}}\to \underline{\Spec}_X\left(\bigoplus_{\underline{n}\in\Z^d}L_{(M_1n_1,\, ...,\,M_dn_d)}\right)\]
where $\mathbb{G}_m^d$ acts on the right hand side with stabilizer being a subgroup of $\prod \mu_{M_i}$. Therefore, the stabilizers on $L_{(M_1n_1,\, ...,\,M_dn_d)}$ are finite subgroups of the torus $\mathbb{G}_m^d$. By definition, when restricting $\Ln$ to the fiber $X_b$ for any $\underline{n}$, they are reflexive. So a fiber in $\mathcal{P}_{\{\Ln\}}$ over $b\in B$ is the spectrum of an algebra with reflexive components over an S2 base, hence it's S2. Since $\mathcal{P}_{\{\Ln\}}\to \mathcal{X}_{\{\Ln\}}$ is smooth and surjective, the quotient stack $\mathcal{X}_{\{\Ln\}}|_{b}$ is also S2.\medskip\\
For (2), note that the invariant part of the $\OO_X$-algebra $\bigoplus_{\underline{n}\in\Z^d}\Ln$ is $L_{\underline{0}}=\OO_X$, so $\mathcal{X}_{\{\Ln\}}\to X$ is the coarse moduli space. For (3), note that $\mathcal{L}_1^{n_1}\otimes\cdots\otimes \mathcal{L}_d^{n_d}$ is the degree $\underline{n}$ component of the $\OO_{\mathcal{X}}$-algebra $\bigoplus_{\underline{n}\in\Z^d}\mathcal{L}_1^{n_1}\otimes\cdots\otimes \mathcal{L}_d^{n_d}$ whose spectrum over $\mathcal{X} $ is $\mathcal{P}_{\{\Ln\}}$. $\Ln$ is the degree $\underline{n}$ component of $\bigoplus_{\underline{n}\in\Z^d}\Ln$, whose spectrum over $X$ is $\mathcal{P}_{\{\Ln\}}.$\medskip\\
For (4), by construction we have\begin{align*}
    \qquad\qquad\mathcal{P}_{\{\Ln'\}} & = \underline{\Spec}_{X'}\bigoplus_{\underline{n}\in\Z^d} \Ln'\\
    & \cong \underline{\Spec}_{X'}\varphi^*\bigoplus_{\underline{n}\in\Z^d} \Ln \qquad \text{($\alpha_{\underline{n}}$ are isomorphisms compatible with multiplication)}\\
    & \cong \varphi^*\underline{\Spec}_{X}\bigoplus_{\underline{n}\in\Z^d} \Ln = \varphi^*\mathcal{P}_{\{\Ln\}}
    \end{align*}
    also $\mathcal{X}_{\{\Ln'\}}\cong\varphi^*\mathcal{X}_{\{\Ln\}}$ because $\varphi^*$ respects the $\Z^d$-grading so we can take $\mathbb{G}_m^d$-quotient.\qed
\end{pf}
\begin{thm}
    The category of Koll\'ar families of Seifert $\mathbb{G}_m^d$-bundles is equivalent to the category of $d$-uniformized $d$-cyclotomic orbispace via the base preserving functors
    \[(X\to B,\{\Ln\}_{\underline{n}\in\Z^d})\mapsto (\mathcal{X}_{\{\Ln\}},\mathcal{L}_1,\,...,\,\mathcal{L}_d)\] where $\mathcal{X}_{\{\Ln\}}=[\mathcal{P}_{\{\Ln\}} /\mathbb{G}_m^d]$, and its inverse is given by the coase moduli map
    \[(\mathcal{X}\to B,\mathcal{L}_1,\,...,\, \mathcal{L}_d)\mapsto \left(X\to B, \Ln=(\pi_*\mathcal{L}_1^{\otimes n_1})[\otimes]\cdots[\otimes](\pi_*\mathcal{L}_d^{\otimes n_d})\right)\]
    \label{seiferteq}
\end{thm}
\begin{rem}
    The difference between the Koll\'ar category of $d$ $\Q$-line bundle sequences and the Koll\'ar category of $d$ $\Q$-line bundles is precisely given by the one extra condition on the stack level: if we insist that the coarse moduli map is an isomorphism away from a codimension $2$ locus, then one can show that all $L_{n_i\underline{e}_i}$ are the $n_i$-th reflexive power of $L_{\underline{e}_i}$. In Section 3, we called a family satisfying this condition a family of $d$-uniformized twisted variety.
\end{rem}
\begin{pf}
    The above proposition gives a functor from the Koll\'ar category of Seifert $\mathbb{G}_m^d$-bundles to the category of $d$-uniformized $d$-cyclotomic orbispaces. Conversely, we take a $d$-uniformized $d$-cyclotomic orbispace $(\mathcal{X}\to B,\mathcal{L}_1,\,...,\, \mathcal{L}_d)$, with coarse moduli space $\pi:\mathcal{X}\to X$. Since $\mathcal{X}$ is a tame stack, the formation of coarse moduli space commutes with arbitrary base change; by universal property of the coarse moduli space, all fibers of $X\to B$ are reduced. If $M_i$ is the index of $\mathcal{I}_{\mathcal{X}}\hookrightarrow T\times \mathcal{X}$ with respect to the $i$-th component of $T$, then $\mathcal{L}_i^{M_i}$ descends to a line bundle on the coarse space $X$, this is exactly $\pi_*\mathcal{L}_i^{M_i}$.\medskip\\
    Around any geometric point of $\mathcal{X}$, there's an \'etale neighbourhood isomorphic to the stack quotient $[V/G]$ where $V$ is an S2 affine scheme and $G$ is a product of finite $\mu_{m_i}$ groups. Locally, the coarse moduli space is just $V/G$, which is still S2 since the invariants are direct summands in $\OO_V$. Recall that the coarse moduli map of a tame stack is good. Hence for any base change to $T$,
    \[\begin{tikzcd}
	{\mathcal{X}_T} & {\mathcal{X}} \\
	{X_T} & X \\
	T & B
	\arrow["{\widetilde{\varphi}}", from=1-1, to=1-2]
	\arrow["p", from=1-1, to=2-1]
	\arrow["{\text{CMS}}"', from=1-1, to=2-1]
	\arrow["\lrcorner"{anchor=center, pos=0.125}, draw=none, from=1-1, to=2-2]
	\arrow["\pi"', from=1-2, to=2-2]
	\arrow["{\text{CMS}}", from=1-2, to=2-2]
	\arrow["\varphi", from=2-1, to=2-2]
	\arrow[from=2-1, to=3-1]
	\arrow["\lrcorner"{anchor=center, pos=0.125}, draw=none, from=2-1, to=3-2]
	\arrow[from=2-2, to=3-2]
	\arrow["g", from=3-1, to=3-2]
\end{tikzcd}\]
    by theorem 4.7 of \cite{gms}, for any line bundle $\mathcal{L}$ on $\mathcal{X}$, we have
    \[\varphi^*\pi_*\mathcal{L}=p_*\widetilde{\varphi}^*\mathcal{L}.\]
    This implies that the restriction to fibers is the same as the fiberwise construction. The sheaves $\pi_{b,*}(\prod_{i=1}^d\mathcal{L}_i^{n_i}|_b)$ are direct summands of the algebra of $\mathcal{P}_{\{\Ln|_b\}}$, which is affine over $X_b$ with S2 fibers. So $\Ln|_b$ is S2 on the reduced scheme $X_b$ of pure dimension, 5.1.1 of \cite{AH} implies that the $\Ln|_b$'s are reflexive. It's clear that $\Ln\otimes L_{M_i\underline{e}_i}\cong L_{\underline{n}+M_i\underline{e}_i}$. Indeed, by the projection formula, it suffices to show that the natural map $\pi^*\pi_*\mathcal{L}_i^{M_i}\to \mathcal{L}_i^{M_i}$ is an isomorphism. \'Etale locally, we have the chart for the coarse space map
   \[\begin{tikzcd}
	{[\Spec A/G]} & {\mathcal{X}} \\
	{\Spec A^G} & X
	\arrow["{\text{\'et}}", from=1-1, to=1-2]
	\arrow[from=1-1, to=2-1]
	\arrow["\lrcorner"{anchor=center, pos=0.125}, draw=none, from=1-1, to=2-2]
	\arrow["\pi", from=1-2, to=2-2]
	\arrow["{\text{\'et}}"', from=2-1, to=2-2]
\end{tikzcd}\]
Now if $N$ is a locally free $A$-module, then for any prime $\mathfrak{p}\subset A$, take $\mathfrak{q}=A^G\cap \mathfrak{p}$, which is a prime ideal in $A^G$. We have $$(N^G\otimes_{A^G}A)_{\mathfrak{p}}=N^G_{\mathfrak{q}}\otimes_{A^G_{\mathfrak{q}}}A_{\mathfrak{p}}\cong A^G_{\mathfrak{q}}\otimes_{A^G_{\mathfrak{q}}}A_{\mathfrak{p}}\cong A_{\mathfrak{p}}$$
since localization respects tensor, preserves exactness and $N^G_{\mathfrak{q}}$ is free. $\mathcal{L}_i^{M_i}$ descends to $X$, so the actions on both sides are trivial. Thus, the map $\pi^*\pi_*\mathcal{L}_i^{M_i}\to \mathcal{L}_i^{M_i}$ is locally isomorphic, and hence isomorphic.
Assume that $U$ is the finite intersection of all opens $U_{\underline{n}}$ where $\Ln$ is invertible for $0\leq n_i<M_i-1$ for every $i$. Then all $\Ln$'s are invertible over $U$.\medskip\\
    Finally, regarding the multiplication structure, denote by $\mathcal{J}_1,...,\mathcal{J}_d$ the uniformizing line bundles on $\mathcal{X}_T$. Clearly $\mathcal{J}_i=\widetilde{\varphi}^*\mathcal{L}_i$ for each $i$. Then,
    \[p_*\mathcal{A}^{\underline{n}}=p_*\widetilde{\varphi}^*\mathcal{L}^{\underline{n}}=\varphi^*\pi_*\mathcal{L}^{\underline{n}}\]
    The multiplication maps $p_*\mathcal{A}^{\underline{n}}\otimes p_*\mathcal{A}^{\underline{m}}\to p_*\mathcal{A}^{\underline{n}+\underline{m}}$ are just the pullbacks of the multiplication maps on $X$ $\varphi^*(\Ln\otimes\Lm)\to\varphi^*\Lmn$.\qed
\end{pf}
\subsection{Moduli space of Seifert $\mathbb{G}_m^d$-bundles} Due to the main theorem, we have
\[\begin{tikzcd}
	{\{\text{Seifert }\mathbb{G}_m^d\text{-bundles}\}} & {\{d\text{-uniformized }d\text{-cyclotomic orbispaces}\}}
	\arrow["{1:1}", tail reversed, from=1-1, to=1-2]
\end{tikzcd}\]
When we fix the Hilbert function $\mathfrak{F}$, then by 4.1.2 of \cite{AH}, 
$$\text{Orb}_{\mathfrak{F}}^{\mathcal{L}_,\,...,\,\mathcal{L}_d,\,\underline{\alpha}}\subset \text{Sta}_{\mathfrak{F}}^{\mathcal{L}_,\,...,\,\mathcal{L}_d,\,\underline{\alpha}}$$
is an open substack. Note that this is exactly the moduli of the left hand side of our correspondence:\begin{prop}
    The moduli stack of all Seifert $\mathbb{G}_m^d$-bundles $(X\to B,\{\Ln\})$ with $\chi(X_b,\Ln)=\mathfrak{F}(\underline{n})$ for any $b\in B$ and for all $\underline{n}\in\Z^d$ is algebraic and locally of finite type. \label{modulisft}
\end{prop}
\noindent One may also wish to fix the coarse space family $X\to B$ and study the Seifert $\mathbb{G}_m$-bundles on it. For this purpose, observe that $X\to B$ admits a relatively very ample line bundle $L_1^{\alpha_1}\otimes\cdots\otimes L_d^{\alpha_d}$, embedding $X$ into $\PP^n_B$, and each fiber $X_b$ has Hilbert function $\chi(X_b,L_{s\underline{\alpha}})=\mathfrak{F}(s\underline{\alpha})$. This gives a map
$$\text{Orb}_{\mathfrak{F}}^{\mathcal{L}_,\,...,\,\mathcal{L}_d,\,\underline{\alpha}}\to \left[\hilb(\PP^n,\mathfrak{F(s)})\middle/ \text{PGL}_{n+1}\right]$$
So the moduli of Seifert $\mathbb{G}_m^d$-bundles on a fixed family $X\to B$ is simply taking a fiber of this map. This shows that:
\begin{prop}
    The moduli space of all Seifert $\mathbb{G}_m^d$-bundles $\{\Ln\}$ on $X\to B$ with $\chi(X_b,\Ln)=\mathfrak{F}(\underline{n})$ for any $b\in B$ and for all $\underline{n}\in\Z^d$ is algebraic and locally of finite type.
\end{prop}

\subsection{Special case where all fibers $X_b$ are normal}
If we assume that all the fibers are normal, then we can obtain a few nice properties. In particular, the construction $\mathcal{P}_{\{\Ln\}}\to X$ satisfies the Definition 1 of \cite{Seifert} and hence it is a ``Seifert $\mathbb{G}_m$-bundle'' when $d=1$. Conversely, we have\begin{thm}
    [c.f. Theorem 7.1 of \cite{Seifert}] If $X$ is a normal variety and $f:Y\to X$ is a Seifert $\mathbb{G}_m$-bundle, then for unique choice of a rank 1 reflexive sheaf $L$ and distinct irreducible divisors $D_i$, $s_i\in (0,1)\cap\Q$, we have
    \begin{align*}
        L_n &\coloneqq L^{[n]}\left(\sum_i\left\lfloor  ns_i \right\rfloor D_i\right) \qquad\text{and}\\
        Y &\cong \underline{\Spec}_X\left(\oplus_{n\in\Z}L_n\right)
    \end{align*}
    where the multiplication structure is given by the obvious map 
    \[L^{[n]}\left(\sum_i \lfloor ns_i\rfloor D_i\right)\otimes L^{[m]}\left(\sum_i \lfloor ms_i\rfloor D_i\right)\to L^{[m+n]}\left(\sum_i \lfloor (m+n)s_i\rfloor D_i\right)\]
    which makes sense since $\lfloor ms_i\rfloor+\lfloor ns_i\rfloor\leq \lfloor (m+n)s_i\rfloor$.
\end{thm}
\noindent This implies that under the normality assumption, the construction $\mathcal{P}_{\{\Ln\}}\to X$ coincides with the notion of Seifert $\mathbb{G}_m$-bundles. Parts (2) and (3) of Theorem 7 in \cite{Seifert} obviously follow from our discussion above. In fact, this 1-to-1 correspondence holds in more generality due to Definition 9.50 and Theorem 9.51 of \cite{singmmp}. $X$ may be assumed to be only a seminormal algebraic space, and the sheaves $L_i$ are just torsion-free sheaves on $X$ of rank 1 or 0 at the generic points. Yet the Seifert $\mathbb{G}_m$-bundles $Y\to X$ ($Y$ is reduced) is the same information as $(X,\{L_i\})$. \medskip\\
\noindent When all fibers $X_b$ in the family $X\to B$ are normal, we may take any codimension 1 point $q$ in $X_b$, then $(X_b)_q$ is the spectrum of a DVR $R$, with uniformizer $t$ because $X_b$ is regular in codimension 1. 
\[\begin{tikzcd}
	{(\mathcal{X}_b)_q}=[\Spec A/G] \\
	{(X_b)_q=\Spec R}
	\arrow["{\text{CMS}}", from=1-1, to=2-1]
\end{tikzcd}\]
Here $\Spec A$ is reduced by 9.51 of \cite{singmmp}, and this coarse space map is a finite morphism. Hence every associated prime $\mathfrak{p}$ of $A$ must be sent to the general point of $\Spec R$. In particular, every line bundle on the source, as a torsion-free sheaf on $\Spec A$, must also be torsion-free over $R$. Since $R$ is a DVR, being torsion-free is the same as being a line bundle.
\begin{rem} If we further impose the condition that the stack fibers $\mathcal{X}_b$ are all normal (which always cuts out an open locus in the moduli), then $(\mathcal{X}_b)_q$ would also be smooth. By definition of orbispaces, this map is only ramified at the Cartier divisor $(t)$ (the closed point of DVR). Theorem 1 of \cite{bottomup} implies that the source can be more concretely described as the root stack $$\sqrt{\Spec R,(t)}=\left[\Spec \frac{R[s]}{(t=s^n)}\middle/\mu_n\right]$$
\end{rem} 
\noindent As we have already argued, for any $\mathcal{L}^{\underline{n}}$, there exists an open subset $U_{\underline{n}}\subset X$, relatively big over $B$, where the pushforward $\pi_* \mathcal{L}^{\underline{n}}=\Ln$ is invertible on $U$. We may take the intersection $U$ of all $U_{\underline{n}}$ for $0\leq n_i<M_i-1$ for every $i$. This $U$ would remain relatively big, and satisfy the condition that for all $\underline{n}$, $\Ln$ is invertible on $U$, yielding:
\begin{prop}
    If $(f:X\to B,\{\Ln\})$ is a Koll\'ar family of Seifert $\mathbb{G}_m^d$-bundles, and all fibers $X_b$ are normal, then there exists a relatively big open $U$ where all the $\Ln$'s are invertible.
\end{prop}
\begin{prop}
    Under the assumption in the previous proposition, all $\Ln$'s are reflexive as coherent sheaves on $X$.
\end{prop}
\begin{pf}
    We use the condition (1) of the definition of the Koll\'ar family of Seifert $\mathbb{G}_m^d$-bundles and the above proposition. In fact, the restriction of any $\Ln$ to any fiber is reflexive and hence S2. By Corollary 3.8 in \cite{HassettKovacs2004}, local freeness of $\Ln$ on $U$ implies that they are reflexive on $X$. \qed
\end{pf}
\begin{cor}
    All the multiplication maps must factor through the reflexive product uniquely.
 \[\begin{tikzcd}
	{\Ln \otimes\Lm} && \Lmn \\
	& {\Ln [\otimes]\Lm}
	\arrow[from=1-1, to=1-3]
	\arrow[from=1-1, to=2-2]
	\arrow[from=2-2, to=1-3]
\end{tikzcd}\]
\end{cor} 
\bibliographystyle{alpha}
\bibliography{ref.bib}

\end{document}